\documentclass[a4paper,11pt]{amsart}
\usepackage{amssymb,amsfonts,amsmath}
\def\d{\delta}

\def\C{\mathbb{C}}
\def\c2{\mathbb{C}^2}
\def\R{\mathbb{R}}

\def\Z{\mathbb{Z}}
\def\N{\mathbb{N}}
\def\P{\mathbb{P}}

\def\1{\bold{1}}

\def\a{\alpha}
\def\b{\beta}
\def\e{\varepsilon}
\def\l{\lambda}
\def\f{\varphi}
\def\g{\gamma}

\def\p{\psi}

\def\o{\omega}

\def\D{\overline{\partial}}

\newtheorem{lem}{Lemma}[section]
\newtheorem{pro}[lem]{Proposition}
\newtheorem{defi}[lem]{Definition}
\newtheorem{def/not}[lem]{Definition/Notations}

\newtheorem{thm}[lem]{Theorem}
\newtheorem{cor}[lem]{Corollary}
\newtheorem{rqe}[lem]{Remark}

\newtheorem{exa}[lem]{Example}
\newtheorem{exas}[lem]{Examples}

\newenvironment{proof3.1}
{\noindent {\it{Proof of theorem 3.1}}}{$\Box$ \linebreak[4]}
\newenvironment{sketch}
{\noindent {\it{Sketch of proof.}}}{$\Box$ \linebreak[4]}

\begin{document}

\title[Monge-Amp\`ere operators on compact K\"ahler manifolds]
{Monge-Amp\`ere operators on compact K\"ahler manifolds}

\author{Vincent GUEDJ \& Ahmed ZERIAHI}

\begin{abstract}
We study the complex Monge-Amp\`ere operator on compact K\"ahler manifolds.
We give a complete description of its range on the set of $\o$-psh functions
with $L^2$-gradient and finite self-energy, generalizing to this compact setting
results of U.Cegrell from the local pluripotential theory.
We give some applications to complex dynamics and to the existence of K\"ahler-Einstein
metrics on singular manifolds.
\end{abstract}

\maketitle


\section*{Introduction}

 Let $X$ be a compact connected K\"ahler manifold, $\dim_{\C} X=k$, equipped
 with a K\"ahler form $\o$. Given a positive Radon measure $\mu$ on $X$, we
 want to study the complex Monge-Amp\`ere equation
$$
\text{ }
\!\!\!\!\!\!  \!\!\! \!\!\! \!\!\! \!\!\! \!\!\! 
[MA](X,\o,\mu) \hskip2cm
(\o+dd^c \f)^k=\mu    \;  ,
$$
where $\f$, the unknown function, is such that $\o_{\f}:=\o+dd^c \f$ is a positive current.
Such functions are called $\o$-plurisubharmonic. We refer the reader to
[16] for their basic properties. Here $d=\partial+\D$ and $d^c=\frac{1}{2i\pi}(\D-\partial)$.

An obvious necessary condition to solve $[MA]$ is that
$\mu(X)=Vol_{\o}(X):=\int_X \o^k$. In the sequel we always assume $\o$ has been normalized
so that $Vol_{\o}(X)=1$, hence we only consider probability measures.
Trying to solve $[MA]$ one immediately faces two problems. The Monge-Amp\`ere operator
$\o_{\f}^k$ is not well defined on the set $PSH(X,\o)$ of all $\o$-psh functions, hence
one has to restrict to subclasses of $\o$-psh functions.
In the whole article we only consider $\o$-psh functions with $L^2$-gradient: this is the
class ${\mathcal E}(X,\o)$, on which the Monge-Amp\`ere operator is well-defined
when $k=\dim_{\C} X=2$,
as was already observed by E.Bedford and A.Taylor in [4]
(see [7] for a more recent account). 

The second problem is that solutions to $[MA]$ are far from being unique, e.g. 
if $\mu$ charges points. One has to further restrict to subclasses of ${\mathcal E}(X,\o)$
in order to give an interesting description of the range of the Monge-Amp\`ere operator.
The first and cornerstone  result in this direction is due to S.-T.Yau [25] who proved 
(see also [1]) the following

\begin{thm}[Yau 1978]
If $\mu$ is a smooth volume form, then there exists a unique 
$\f \in PSH(X,\o) \cap {\mathcal C}^{\infty}(X)$ such that
$$
\o_{\f}^k=\mu \; \; \text{ and } \; \; \sup_X \f=-1.
$$
\end{thm}

This is the solution to a celebrated conjecture of E.Calabi [8], and it has important consequences
in differential geometry (see [22], [23]).

>From the point of view both of complex geometry and pluripotential theory, it is important
to solve $[MA]$ for singular measures $\mu$. A major contribution was made by S.Kolodziej
[18], [19], who proved the following result.

\begin{thm}[Kolodziej, 1998]
If $\mu=f \o^k$ has density $f \in L^p(X), \; p>1$, then there exists a unique 
$\f \in PSH(X,\o) \cap {\mathcal C}^0(X)$ such that
$$
\o_{\f}^k=\mu \; \; \text{ and } \; \; \sup_X \f=-1.
$$
\end{thm}

For applications in complex geometry and complex dynamics one needs to allow unbounded solutions 
$\f$ to Monge-Amp\`ere equations $[MA]$ with more singular measures $\mu$
(see [12], [13] and [2], [15]). This is the main goal of this article. We consider
the subclass ${\mathcal E}^1(X,\o)$ of $\o$-psh functions $\f$ 
such that $\o_{\f}^k$ is well defined and for which
$\f \in L^1(\o_{\f}^k)$. Our main result gives a full
characterization of the Monge-Amp\`ere operator on this class.

\begin{thm}[MAIN THEOREM]
There exists a unique $\f \in {\mathcal E}^1(X,\o)$ such that
$$
\o_{\f}^k=\mu \; \; \text{ and } \; \; \sup_X \f=-1.
$$
if and only if ${\mathcal E}^1(X,\o) \subset L^1(\mu)$.
\end{thm}

The class ${\mathcal E}^1(X,\o)$ contains many unbounded functions, however
these are not too singular (e.g. they have zero Lelong numbers), as follows from 
the condition $\f \in L^1(\o_{\f}^k)$. A similar result was proved by
U.Cegrell [9] in a local setting (for bounded hyperconvex domains in $\C^n$).
Our work originated from an attempt to understand Cegrell classes in the global
context of $\o$-psh functions on compact K\"ahler manifolds. We would like to emphasize
that surprinsingly the global and the local theory are  quite different:

-There is no boundary condition in the compact setting. Integration by parts are
much simpler, as well as compactness arguments. In particular Monge-Amp\`ere masses
are uniformly controlled by $Vol_{\o}(X)=1$. 

-The counterpart is that one looses the homogeneity of the Monge-Amp\`ere operator
$\o_{\f}^2$. If $\f \in PSH(X,\o)$ then $\l \f \in PSH(X,\o)$ for $0 \leq \l \leq 1$
but not for $\l>1$, and one has to get control on  mixed terms 
$\o^2,\o \wedge \o_{\f}, \o_{\f}^2$, due to the affine part of $\o_{\f}=\o+dd^c \f$.
This is an important source of difficulty.

\vskip.2cm

In the major part of the article we 
restrict ourselves to the case of complex surfaces ($k=2$) because
it greatly simplifies both the definition of the Monge-Amp\`ere operator
and the computations. 
However most of our results hold on any 
$k$-dimensional compact K\"ahler manifold, as we explain in section 7.1.
We also assume throughout the article that
$\o$ is a {\it Hodge form}, i.e. a K\"ahler form with integer 
cohomology class.
In this case one can easily regularize $\o$-psh functions with no loss
of positivity (see Appendix in [16]). All our results are
true when $\o$ is merely a K\"ahler form, but our estimates then have to be modified
by a uniform constant.

\vskip.2cm

We now describe the contents of the article more precisely.
In {\it section 1} we define and study the class ${\mathcal E}(X,\o)$ of $\o$-psh functions whose
gradient is in $L^2(X)$. We then define the Monge-Amp\`ere operator $\o_{\f}^2$
in {\it section 2} and prove basic continuity results. The class
${\mathcal E}^1(X,\o)$ is introduced in {\it section 3}. It is a starshaped
convex set, stable under taking maximum (proposition 3.2). These properties are established through
integration by parts which are justified thanks to an important continuity result (Theorem 3.1).
We also prove (Theorem 3.4) that solutions to the Monge-Amp\`ere equation
$[MA]$ are unique -up to an additive constant- in the class ${\mathcal E}^1(X,\o)$.

In {\it section 4} we define and study several intermediate classes ${\mathcal E}^p(X,\o)$,
$p \geq 1$. Not only are they interesting in themselves, but we need to solve first
$[MA]$ with solutions in ${\mathcal E}^p(X,\o)$, $p>1$, before producing
solutions in ${\mathcal E}^1(X,\o)$. We prove our main result (Theorem 5.1) in
{\it section 5} where we give a complete characterization of the range of
the Monge-Amp\`ere operator $\o_{\f}^2$ on each class ${\mathcal E}^p(X,\o)$, $p \geq 1$.
In {\it section 6} we give several examples of functions in ${\mathcal E}^p(X,\o)$
and measures of Monge-Amp\`ere type by linking these integrability properties
to the size of the sublevel sets $(\f<-t)$, estimated in terms of the Monge-Amp\`ere
capacity $Cap_{\o}$. Finally in {\it section 7} we describe how to adapt our
arguments to higher dimension and give some applications to complex dynamics
and to the existence of K\"ahler-Einstein metrics on singular manifolds.

\section{The gradient of qpsh functions}

>From now on -- until section 7 --
we assume $X$ is a compact complex projective {\it surface}
(i.e. of complex dimension 2) and $\o$ is a Hodge form on $X$
normalized so that $\text{Vol}_{\o}(X):=\int_X \o^2 =1$. 
Recall that
$$
PSH(X,\o):=\left\{ \f \in L^1(X) \, / \, dd^c \f \geq -\o \text{ and } \f \text{ is u.s.c. } \right\}
$$
is the set of $\o$-psh functions (see [16]). 
We set
$$
{\mathcal E}(X,\o):=\left\{ \f \in PSH(X,\o) \, / \, 
\f \in L^1(\o \wedge \o_{\f}) \right\}
$$
and 
$$
W^{1,2}(X):=\left\{ \f \in L^2(X) \, / \, \nabla \f \in L^2(X) \right\},
$$
endowed with the Sobolev norm
$||\f||_{W^{1,2}}=||\f||_{L^2}+||\nabla \f||_{L^2}$. 
To simplify notations all $L^p$-norms are computed with respect to the volume
form $\o^2$ unless otherwise specified.
Here the $L^2$-norm of the gradient simply means
$$
||\nabla \f ||_{L^2}:=\left( \int_X d\f \wedge d^c \f \wedge \o \right)^{1/2}.
$$
We shall denote ${\mathcal E}(X,\o)$ by ${\mathcal E}$ when no confusion
can arise.
Recall that every $\o$-psh function can be approximated by a decreasing sequence
of smooth $\o$-psh functions (see Appendix in [16]). This motivates the following

\begin{lem}
Let $\f_j,\f \in  {\mathcal E}(X,\o)$.

1) If $\f_j$ decreases towards $\f$, then $\f_j \rightarrow \f$ in the
Sobolev $W^{1,2}$-norm.

2) If $\f_j \rightarrow \f$ in the Sobolev $W^{1,2}$-norm, then
$\f_j \o_{\f_j} \rightarrow \f \o_{\f}$ and  
$d\f_j \wedge d^c \f_j \rightarrow d\f \wedge d^c \f$ 
in the sense of currents.

\end{lem}

\begin{proof}
Assume first that $\f_j \rightarrow \f$ in the Sobolev norm.
Then $d \f_j \wedge d^c \f_j \rightarrow d\f \wedge d^c \f$ 
and $\f_j d^c \f_j \rightarrow \f d^c \f$ in $L^1$ hence in
the sense of currents. Let $\theta$ be a smooth test form. We get
$$
<\f_j dd^c \f_j, \theta>=<d \f_j \wedge d^c \f_j, \theta>
-<\f_jd^c \f_j,d\theta>
\rightarrow <\f dd^c \f, \theta>,
$$
thus $\f_j \o_{\f_j} \rightarrow \f \o_{\f}$ in the sense of currents.

Assume now $\f_j$ decreases towards $\f$. By the monotone convergence theorem,
$\f_j \rightarrow \f$ in $L^2$, so the question is whether
$\nabla(\f_j-\f) \rightarrow 0$ in $L^2$. We have
$$
0 \leq \int_X d(\f_j-\f) \wedge d^c(\f_j-\f) \wedge \o=
\int_X -(\f_j-\f) (\o_{\f_j}-\o_{\f}) \wedge \o \leq
\int_X (\f_j-\f) \o_{\f} \wedge \o.
$$
Now $\f_j-\f \in L^1(\o_{\f} \wedge \o)$, so it follows
from the monotone convergence theorem again that
$\int_X (\f_j-\f) \o_{\f} \wedge \o \rightarrow 0$.
\end{proof}

In the sequel we shall implicitly make computations on smooth approximants and then 
pass to the limit by using lemma 1.1.

\begin{pro} One has
$$
{\mathcal E}(X,\o)=PSH(X,\o) \cap W^{1,2}(X).
$$ 
Moreover ${\mathcal E}(X,\o)$ is a closed subspace
of $W^{1,2}(X)$.
\end{pro}

\begin{proof}
Let $\f \in PSH(X,\o)$. We can assume w.l.o.g. $\f \leq 0$. 
Let us recall that a $\o$-psh function is $L^p$-integrable for all $p \geq 1$
and has gradient in $L^p$ for all $p<2$ (see [17] and inequality (2)
below). It follows therefore from H\"older inequality that 
$\f d^c \f$ is a well defined current of degree 1 with $L^1$ coefficients.
Observe that
$$
d\left(\f d^c\f \right)+\f \o=d\f \wedge d^c \f+ \f \o_{\f},
$$
where one of the currents on the right hand-side is well defined
as soon as the other is.
It follows therefore from Stokes
theorem that
$$
0 \leq \int_X (-\f) \o_{\f} \wedge \o =||\f||_{L^1}+
\int_X d\f \wedge d^c \f \wedge \o.
$$
Thus $\f \in {\mathcal E}(X,\o)$ iff $\nabla \f \in L^2(X)$.

Clearly if $\f_j \in {\mathcal E}(X,\o)$ converges in the Sobolev
$W^{1,2}$-norm towards $\f$ then $\f_j \rightarrow \f$ in $L^2$ so that
$\f \in PSH(X,\o)$ and $\nabla \f \in L^2(X)$ so that 
$\f \in {\mathcal E}(X,\o)$, i.e. ${\mathcal E}(X,\o)$ is closed in
$W^{1,2}(X)$.
\end{proof}

\begin{pro} \text{ }

1) If $\f,\p \in {\mathcal E}(X,\o)$ then 
$\p \in  L^1(\o_{\f} \wedge \o)$. Therefore $\o_{\f} \wedge \o_{\p}$ is a well defined
probability measure.

2) The set ${\mathcal E}(X,\o)$ is  star-shaped and convex.

3) Assume $\f,\p \in PSH(X,\o)$ with $\f \leq \p$. If 
$\f \in  {\mathcal E}(X,\o)$
then $\p \in {\mathcal E}(X,\o)$. In particular ${\mathcal E}(X,\o)$ is stable
under taking maximum.
\end{pro}

\begin{proof}
Fix $\f, \p \in {\mathcal E}$. We can assume $\p \leq 0$. By Stokes
theorem and Cauchy-Schwarz inequality, we get
\begin{eqnarray*}
\lefteqn{0 \leq \int (-\p) \o_{\f} \wedge \o
=\int (-\p) \o^2 +\int (-\p) dd^c \f \wedge \o} \\
&& =\int (-\p) \o^2+ \int d\p \wedge d^c \f \wedge \o \\
&& \leq \int (-\p) \o^2 +\left( \int d\p \wedge d^c \p \wedge \o \right)^{1/2}
\cdot \left( \int d\f \wedge d^c \f \wedge \o \right)^{1/2} <+\infty.
\end{eqnarray*}
The current $\p \o_{\f}$ is therefore well defined, hence so is
$$
\o_{\f} \wedge \o_{\p}:=\o \wedge \o_{\f} +dd^c \left( \p \o_{\f} \right).
$$
This yields a probability measure, as can be seen by approximating $\f$ and $\p$ by smooth
approximants and by using lemma 1.1.

We now show that ${\mathcal E}$ is convex.
Given $\f,\p \in {\mathcal E}$, it suffices to check that
$u=(\f+\p)/2$ also belongs to ${\mathcal E}$. By symmetry we only need
to show that $\f \in L^1(\o_u \wedge \o)$. Since $\o_u=(\o_{\f}+\o_{\p})/2$
and $\f \in L^1(\o_{\f} \wedge \o)$, this boils down to check that
$\f \in L^1(\o_{\p} \wedge \o)$, which is nothing but 1.3.1.

Assume $\f \in {\mathcal E}$ and $\l \in [0,1]$. Then 
$\l \f \in PSH(X,\o)$ since
$\o_{\l \f}=\l \o_{\f}+(1-\l) \o \geq 0$. Also
$\nabla (\l \f)=\l \nabla\f \in L^2(X)$, hence 
$\l \f \in {\mathcal E}$, i.e. ${\mathcal E}$ is star-shaped.

It remains to prove 3).
Let $\f,\p \in PSH(X,\o)$ with $\f \leq \p$. We assume again
$\p \leq 0$. It follows then from Stokes theorem that
\begin{eqnarray*}
0 \leq \int_X (-\p) \o_{\p} \wedge \o & \leq &
\int_X (-\f) \o_{\p} \wedge \o=\int_X (-\f) \o^2
+\int_X (-\f)dd^c \p \wedge \o \\
&=& \int_X (\p-\f) \o^2+\int_X (-\p) \o_{\f} \wedge \o \\
&\leq& ||\p-\f||_{L^1}+\int_X (-\f) \o_{\f} \wedge \o.
\end{eqnarray*}
Therefore $\f \in {\mathcal E} \Rightarrow \p \in {\mathcal E}$.
\end{proof}

\begin{exas}
\text{ }

1) $L^{\infty}(X) \cap PSH(X,\o) \subset {\mathcal E}(X,\o)$.

\noindent Indeed ${\mathcal E}$ obviously contains constant functions
and if $\f \in PSH(X,\o)$ is bounded from below by some
constant $C$ then $\f \in {\mathcal E}$ by the previous proposition.
Alternatively it may be useful to note the following inequality: if
$0 \leq \f \leq 1/2$, $\f \in PSH(X,\o)$, then 
$$
dd^c (\f^2)=2 d\f \wedge d^c \f+2 \f dd^c \f \geq -\o
$$
hence $\f^2 \in PSH(X,\o)$ with
$$
||\nabla \f||_{L^2}^2 = \int_X d\f \wedge d^c \f \wedge \o
\leq \frac{1}{2} \int_X \o_{\f^2} \wedge \o=\frac{1}{2}
$$

More generally if $\chi:\R \rightarrow \R$
satisfies $\chi''> 0$ and $0 \leq \chi' \leq A$ on $\f(X)$, then
$\chi \circ \f \in {\mathcal E}(X,A \o)$ since
\begin{equation}
dd^c(\chi \circ \f)=\chi'' \circ \f d\f \wedge d^c \f+\chi' \circ \f dd^c \f.
\end{equation}

2) If $\f \in PSH(X,\o)$ is 
bounded near some ample divisor $D$, 
then $\f \in {\mathcal E}$. 

\noindent Indeed let $\o_D$ be a K\"ahler form cohomologous to $[D]$, the
current of integration along $D$. Integrability against $\o$ is equivalent to
integrability against $\o_D$; for simplicity we assume $\o_D=\o$.
We can find $\o'$ a smooth
positive closed $(1,1)$ form cohomologous to $\o$ such that $\o' \equiv 0$
outside some small neighborhood $V$ of $D$ where $\f$ is bounded.
Fix $\chi \geq 0$ smooth such that $\o=\o'+dd^c \chi$ and
assume w.l.o.g. $\f \leq 0$. Then
\begin{eqnarray*}
0 \leq \int_X (-\f) \o_{\f} \wedge \o &=&
\int_X (-\f) \o_{\f} \wedge \o'+ \int_X (-\f) \o_{\f} \wedge dd^c \chi \\
&\leq& ||\f||_{L^{\infty}(V)} \int_X \o_{\f} \wedge \o'
+\int_X \chi \o_{\f} \wedge (-dd^c \f ) \\
&\leq& ||\f||_{L^{\infty}(V)}+||\chi||_{L^{\infty}(X)}<+\infty,
\end{eqnarray*}
since $-dd^c \f \leq \o$, $\chi \o_{\f} \geq 0$ and
$\int_X \o_{\f} \wedge \o'=\int_X \o_{\f} \wedge \o=\int_X \o^2=1$.

These examples are analogous to the psh functions
with compact singularities introduced and studied by N.Sibony [21] in
the local theory (see the survey article [11] and references therein). 

3) If $\o_{\f}$ is the current of integration along some complex
hypersurface $H$ of $X$ then $\f \equiv-\infty$ on $H$ hence
$\f \notin {\mathcal E}$. 
One can also produce examples of functions
$\f \in PSH(X,\o) \setminus {\mathcal E}$ which have zero Lelong number
at all points:
let $\f[z_0:z_1:z_2]=\log|z_0|-\log||(z_0,z_1,z_2)||-1 \leq -1$. Then
$\f \in PSH(\P^2,\o)$, where $\o$ denotes the Fubini-Study K\"ahler form on
$X=\P^2$. If $0 \leq \a \leq 1$ then 
$-(-\f)^{\a} \in PSH(\P^2,\o)$ and a straightforward computation shows that
$-(-\f)^{\a} \in {\mathcal E}(\P^2)$ iff $0 \leq \a <1/2$.

\end{exas}

One can generalize example 1.4.3 as the following proposition shows:

\begin{pro}
Assume $\f \in PSH(X,\o)$, $\f \leq -1$ and fix $0 \leq \a <1/2$.

Then $-(-\f)^{\a} \in {\mathcal E}(X,\o)$. 
In particular every locally pluripolar set is included in the $-\infty$ locus
of a function in ${\mathcal E}(X,\o)$.
\end{pro}

\begin{proof}
Set $\f_{\a}:=-(-\f)^{\a}=\chi \circ \f$, where $\chi(t)=-(-t)^{\a}$.
We assume first $0 \leq \a < 1$.
Observe that $0 \leq \chi' \circ \f \leq \a \leq 1$ and 
$\chi" \circ \f=\a(1-\a) (-\f)^{\a-2} \geq 0$ so that by (1),
$\f_{\a}  \in PSH(X,\o)$ with
\begin{equation}
0 \leq \int_X \frac{1}{(-\f)^{2-\a}} d\f \wedge d^c \f \wedge \o \leq
\frac{1}{1-\a}<+\infty.
\end{equation}
It is well-known that the gradient $\nabla \f$ of a $\o$-psh function $\f$
is in $L^p$ for all $p<2$ (but does not belong to $L^2$ in general, see
example 1.4.3). Inequality (2) shows, more precisely, that $\nabla \f$ belongs
to a weighted version of $L^2$.

Observe now that 
$d \f_{\b} \wedge d^c \f_{\b}=\b^2 (-\f)^{2-2\b} d\f \wedge d^c \f$. If we set
$\b=\a/2<1/2$, it therefore follows from (2) that
$$
\int_X d\f_{\b} \wedge d^c \f_{\b} \wedge \o \leq \frac{\b^2}{1-2\b} <+\infty,
$$
so that $\f_{\b} \in {\mathcal E}(X,\o)$.
\end{proof}

This proposition implies, together with proposition 1.3, the following:

\begin{cor}
Let $\f \in {\mathcal E}(X,\o)$. Then
the measure $\o_{\f} \wedge \o$ does not charge pluripolar sets.
\end{cor}

\section{The complex Monge-Amp\`ere operator}

The complex Monge-Amp\`ere operator $\o_{\f}^2$
can be easily defined for functions $\f \in {\mathcal E}(X,\o)$. This was
already observed by Bedford and Taylor  in the local
context (see [4], [5]). 
Indeed if $\f \in {\mathcal E}(X,\o)$ then we set
$$
\o_{\f}^2:=\o \wedge \o_{\f}+dd^c(\f \o_{\f}).
$$
This is a well defined current of maximal bidegree $(2,2)$
which happens to be a probability measure (proposition 1.3.1).
This operator is continuous on decreasing sequences, as follows from lemma 1.1.

In the context of quasiplurisubharmonic functions, convergence in the Sobolev norm  and 
weak convergence are almost the same as the following result shows.

\begin{thm}
Let $\f_j,\f \in {\mathcal E}(X,\o)$. Assume there exists
$\p \in {\mathcal E}(X,\o)$ such that $\f_j \geq \p$ for all $j \in \N$.
Then the following are equivalent:

1) $\f_j \rightarrow \f$ in $W^{1,2}(X)$.

2) $\f_j \rightarrow \f$  and 
$\f_j \o_{\f_j} \rightarrow \f \o_{\f}$ in the sense of distributions.

3) $\f_j \rightarrow \f$ in the sense of distributions and 
$\int \f_j \o_{\f_j} \wedge \o \rightarrow \int \f \o_{\f} \wedge \o$.
\end{thm}

\begin{proof}
Observe that $1) \Rightarrow 2)$ follows from lemma 1.1 and 
$2) \Rightarrow 3)$ is obvious. So it remains to prove $3) \Rightarrow 1)$.
It is a standard consequence of 
quasi-plurisubharmonicity that
$\f_j \rightarrow \f$ weakly iff $\f_j \rightarrow \f$ in $L^2$; moreover
for all $x \in X$, $\limsup \f_j(x) \leq \f(x)$ with equality off a pluripolar
set (see [17]). 
Assume this is the case. It follows that $(\f_j)$ is uniformly bounded 
from above. We can assume $\f_j,\f \leq 0$, hence our assumption yields
a uniform bound
\begin{equation}
\p \leq \f_j \leq \f_j-\f \leq -\f.
\end{equation}
A repeated application of Stokes theorem yields
\begin{eqnarray*}
\lefteqn{ 0 \leq \int_X d(\f-\f_j) \wedge d^c (\f-\f_j) \wedge \o
= \int_X -(\f-\f_j) (\o_{\f}-\o_{\f_j}) \wedge \o} \\
&& =\int_X (\f-\f_j) \o^2 +2 \int_X (\f_j-\f) \o_{\f} \wedge \o
+\int_X (\f \o_{\f}-\f_j \o_{\f_j}) \wedge \o.
\end{eqnarray*}
The first integral converges to $0$ since $\f_j \rightarrow \f$ in $L^2$.
The last one also by assumption 3). Thanks to $(3)$ we can apply Fatou's lemma
and get
$$
\limsup \int_X (\f_j-\f) \o_{\f} \wedge \o \leq 0,
$$
hence actually $\int_X (\f_j-\f) \o_{\f} \wedge \o \rightarrow 0$
and $\f_j \rightarrow \f$ in $W^{1,2}(X)$.
\end{proof}

\begin{cor}
If $\f_j \in {\mathcal E}(X,\o)$ increases towards $\f \in {\mathcal E}(X,\o)$,
then $\f_j \rightarrow \f$ in the Sobolev $W^{1,2}$-norm,
hence $\o_{\f_j}^2 \rightarrow \o_{\f}^2$.
\end{cor}

\begin{proof} 
Assume $(\f_j) \in {\mathcal E}$ is an increasing 
sequence which converges
in $L^2$ towards $\f$. Then 
$\f_j(x) \rightarrow \f(x)$ at every point of $X \setminus P$, where $P$
is a pluripolar set (see [5]). 
We can assume w.l.o.g $\f_j \leq \f \leq 0$.

Fix $p \in \N$ and consider indices $j \geq p$. Then
$$
\int_X (-\f_j) \o_{\f_j} \wedge \o \leq \int_X (-\f_p) \o_{\f_j} \wedge \o
=\int_X (-\f_j) \o_{\f_p} \wedge \o+ \int_X (\f_j-\f_p) \o^2.
$$
Since $\f_j(x) \nearrow \f(x)$ at $\o_{\f_p} \wedge \o$ 
(and $\o^2$) almost every point $x$, we infer
$$
\limsup \int_X (-\f_j) \o_{\f_j} \wedge \o \leq 
\int_X (-\f) \o_{\f_p} \wedge \o+\int_X(\f-\f_p) \o^2=
\int_X (-\f_p) \o_{\f} \wedge \o.
$$
Now $\f_p \nearrow \f$ as $p \rightarrow +\infty$ with pointwise
convergence $\o_{\f} \wedge \o$-almost everywhere. This shows
$$
\limsup_{j \rightarrow +\infty} \int_X (-\f_j) \o_{\f_j} \wedge \o
\leq \int_X (-\f) \o_{\f} \wedge \o.
$$

It follows that the sequence of positive currents $(-\f_j) \o_{\f_j}$ has 
uniformly bounded mass. Let $S$ be a cluster point of this sequence.
We have just shown $||S|| \leq ||(-\f) \o_{\f}||$.
We claim $S \geq (-\f) \o_{\f}$. Indeed if $\theta \geq 0$ is a positive test
form, then
$$
<\f_j \o_{\f_j}, \theta> \, \leq  \, <\f \o_{\f_j},\theta> \,
\leq \, <\f^{\e} \o_{\f_j},\theta>,
$$
where $\f^{\e}$ denotes smooth $\o$-psh functions that decrease towards $\f$.
Since $\o_{\f_j}$ weakly converges towards 
$\o_{\f}$, we infer, 
$$
<S,\theta> \, \geq \, <(-\f^{\e}) \o_{\f},\theta> 
\stackrel{\e \rightarrow 0}{\longrightarrow} <(-\f) \o_{\f} ,\theta>,
$$
therefore $S=\f \o_{\f}$, i.e. $\f_j \o_{\f_j} \rightarrow \f \o_{\f}$.
By theorem 2.1 this shows $\f_j \rightarrow \f$ in the Sobolev 
$W^{1,2}$-norm.
\end{proof}

In [16] we have started to study the   
Monge-Amp\`ere capacity associated to $\o$ 
which is defined as follows:
$$
Cap_{\o}(K):=\sup \left\{ \int_K \o_{u}^2 \; / \, u \in PSH(X,\o),
\, 0 \leq u \leq 1 \right\},
$$
where $K$ is any Borel subset of $X$. This capacity vanishes on pluripolar
sets, more precisely $Cap_{\o}(\f <-t) \leq C_{\f}/t$ for every fixed $\o$-psh
function $\f$. This estimate is sharp in the sense that
$Cap_{\o}(\f <-t) \geq C'_{\f}/t$ when $\o_{\f}$ is the current of integration
along an hypersurface.
However when $\f$ belongs to ${\mathcal E}(X,\o)$, one can establish finer
estimates as the following proposition
shows. 

\begin{pro}
Assume $\f \in {\mathcal E}(X,\o)$. Then there exists $C_{\f}>0$
such that  for all $t>0$,
$$
Cap_{\o}(\f <-t) \leq \frac{C_{\f}}{t^2}.
$$
\end{pro}

\begin{proof}
We can assume w.l.o.g. $\f \leq 0$.
Fix $u \in PSH(X,\o)$ with $-1 \leq u \leq 0$. By Chebyshev inequality,
$$
\int_{(\f<-t)} \o_u^2 \leq \frac{1}{t^2} \int_X \f^2 \o_u^2=
\frac{1}{t^2} \left[ \int_X \f^2 \o^2+ \int_X \f^2 \o \wedge dd^c u
+\int_X \f^2 dd^c u \wedge \o_u \right].
$$
So we need to get an upper bound on the last two integrals that is uniform in
$u$. We first get a bound on $I_u=\int_X \f^2 \o \wedge dd^c u$.
Observe that 
$$
-dd^c (\f^2)=-2 d\f \wedge d^c \f+2(-\f)dd^c \f \leq 2(-\f) \o_{\f}.
$$
It follows therefore from Stokes theorem that 
$$
I_u:=\int_X u dd^c \f^2 \wedge \o \leq 2 \int_X (-u)(-\f)\o_{\f} \wedge \o
\leq 2 \int_X (-\f) \o_{\f} \wedge \o.
$$
Similarly one gets
\begin{eqnarray*}
\!\!\!\!\!\!\!\!\!
\lefteqn{II_u:=\int_X \f^2 dd^c u \wedge \o_u 
\leq   2 \int_X (-\f) \o_{\f} \wedge \o_u} \\
&& \leq   2 \int_X (-\f) \o_{\f} \wedge \o + 2\int_X (-u) \o_{\f}^2 
 \leq  2 \int_X (-\f) \o_{\f} \wedge \o +2.
\end{eqnarray*}
We infer
$$
Cap_{\o}(\f <-t) \leq \frac{1}{t^2}
\left[ \int_X \f^2 \o^2+4 \int_X (-\f) \o_{\f} \wedge \o+2 \right].
$$
\end{proof}

Our aim in this paper is to describe the range of the
Monge-Amp\`ere operator on various subclasses of ${\mathcal E}(X,\o)$.
It is an interesting open question to obtain a description of the set of
probability measures 
$$
{\mathcal M}(X,\o):=\left\{ \o_{\f}^2 \; / \; \f \in {\mathcal E}(X,\o)
\right\}.
$$
One can ask for instance if every probability measure on $X$ belongs to 
${\mathcal M}(X,\o)$ ?
One of the difficulties lies in the lack of uniqueness of solutions
$\f \in {\mathcal E}(X,\o)$ to the equation $\o_{\f}^2=\mu$, as the following
example shows.

\begin{exa}
Let $\o=\o_{FS}$ be the Fubini-Study K\"ahler form on $X=\P^2$.
Let $0$ be the origin in some affine chart $0 \in \C^2 \subset \P^2$.
We use $(z_1,z_2)$ as affine coordinates in $\C^2$ and let
$[z_0:z_1:z_2]$ denote the homogeneous cooordinates on $\P^2$. Consider
$$
\f[z_0:z_1:z_2]:=\int_{\a \in \P^1} \log |\a_1z_1+\a_2z_2| d\nu(\a)
-\log||(z_0,z_1,z_2)||,
$$
where $\nu$ denotes a probability measure on the Riemann sphere
$\P^1=\{\a=[\a_1:\a_2]\}$. Then $\f \in PSH(\P^2,\o)$
and $\f \in {\mathcal E}(\P^2,\o)$ iff $\nu$ has finite self-energy
(i.e. if its potentials have gradient in $L^2$).
Assume this is the case. Observe that $\log|\a_1z_1+\a_2z_2|$ is harmonic on
each radial line through the origin, except at the origin, to conclude that
$$
\o_{\f}^2=\d_0:=\text{Dirac mass at point } 0 \in \C^2 \subset \P^2.
$$
Thus the set of solutions 
$\{\f \in {\mathcal E}(\P^2,\o) \, / \, \o_{\f}^2=\d_0 \}$ has infinite
dimension and contains a subset isomorphic to
$PSH(\P^1,\o_{\P^1}) \cap W^{1,2}(\P^1)$, where $\o_{\P^1}$ denotes here the
Fubini-Study K\"ahler form on the Riemann sphere $\P^1$.
\end{exa}

In the remaining part of this article, we are going to define and study
several subclasses of ${\mathcal E}(X,\o)$ on which solutions of
Monge-Amp\`ere equations are essentially unique. This is the key to the
description of the range of $\o^2_{\f}$ on these classes.

\section{The class ${\mathcal E}^1(X,\o)$}

Our main subject of interest in the sequel is the following class of qpsh functions,
$$
{\mathcal E}^1(X,\o):=\left\{ \f \in {\mathcal E}(X,\o) \, / \, \f \in L^1(\o_{\f}^2) \right\}.
$$
When no confusion can arise, we shall simply denote ${\mathcal E}^1(X,\o)$
by ${\mathcal E}^1$.

Of course bounded $\o$-psh functions belong to ${\mathcal E}^1$. 
We will exhibit in examples 6.3 below 
unbounded functions that belong to ${\mathcal E}^1$. These however have mild singularities: it follows 
from a result of J.-P.Demailly [11] that if $\f \in {\mathcal E}$ has positive Lelong number
at some point $a \in X$, then $\o_{\f}^2$ has some positive Dirac mass at point $a$, so that
$\f$ cannot be integrable with respect to $\o_{\f}^2$.

Before establishing basic properties of the class ${\mathcal E}^1$, we start by proving a useful
continuity result.

\begin{thm} [Continuity]
Let $\f_j \in {\mathcal E}^1(X,\o)$ be a decreasing sequence.
Then the sequence
$(\int_X (-\f_j) \o_{\f_j}^2)_j$ is bounded  if and only if 
$\f:=\lim \searrow \f_j \in {\mathcal E}^1(X,\o)$, and in this case
$$
(-\f_j) \o_{\f_j}^2 \longrightarrow (-\f) \o_{\f}^2 \; \; 
\text{ in the weak sense of currents}.
$$
\end{thm}

\begin{proof}
We can assume w.l.o.g. that $\f_j \leq 0$. Observe that
$$
0 \leq \int_X (-\f_j) \o_{\f_j}^2=\int_X (-\f_j) \o^2+
\int_X d\f_j \wedge d^c \f_j \wedge \o
+\int_X d\f_j \wedge d^c \f_j \wedge \o_{\f_j}.
$$
Since $d\f_j \wedge d^c \f_j \wedge \o_{\f_j} \geq 0$, it follows that
$(\f_j)$ has bounded Sobolev norm, hence
$\f_j \rightarrow \f$ in the Sobolev $W^{1,2}$-norm 
(so $\f \in {\mathcal E}^0(X,\o))$ and $\o_{\f_j}^2 \rightarrow \o_{\f}^2$
by Theorem 2.1.

Assume first that the sequence $(\int_X (-\f_j) \o_{\f_j}^2)$ is bounded.
Let $\nu$ be a cluster point of the sequence of positive measures
$(-\f_j) \o_{\f_j}^2$. It follows from standard arguments that 
$(-\f) \o_{\f}^2 \leq \nu$, hence, in particular, 
$\f \in {\mathcal E}^1(X,\o)$. Indeed we can find smooth $\o$-psh functions
$(\f_j^{\e})_{\e>0}$ that decrease towards $\f_j$ as
$\e \searrow 0$. Let $\chi \geq 0$ be a test
function, then for $j \geq p$, we obtain
$$
<(-\f_j) \o_{\f_j}^2, \chi> \, \geq \, <(-\f_p) \o_{\f_j}^2, \chi>
\, \geq \, <(-\f_p^{\e}) \o_{\f_j}^2, \chi>,
$$
thus
$$
<\nu,\chi> \, \geq \, <(-\f_p^{\e}) \o_{\f}^2,\chi> 
\stackrel{\e \rightarrow 0}{\longrightarrow} <-\f_p \o_{\f}^2,\chi>
\stackrel{p \rightarrow +\infty}{\longrightarrow} <-\f \o_{\f}^2,\chi>,
$$
by the monotone convergence theorem.

We now show that $\nu$ and $(-\f) \o_{\f}^2$ have the same mass. This will
prove that $\nu=(-\f) \o_{\f}^2$ is the unique cluster point, hence
$\f_j \o_{\f_j}^2 \rightarrow \f \o_{\f}^2$. Since $-\f_j \leq -\f$,
it follows from Stokes theorem that
$$
\int_X (-\f_j) \o_{\f_j}^2 \leq \int_X (-\f) \o_{\f_j}^2
=\int_X (\f_j-\f) \o \wedge \o_{\f_j}
+\int_X (-\f_j) \o_{\f} \wedge \o_{\f_j}.
$$
The first integral converges to $0$
since $||\f_j-\f||_{W^{1,2}} \rightarrow 0$.
We estimate the last one by using Stokes theorem again,
$$
\int_X (-\f_j) \o_{\f} \wedge \o_{\f_j} \leq 
\int_X (-\f) \o_{\f} \wedge \o_{\f_j} \leq 
\int_X (\f_j-\f) \o \wedge \o_{\f}
+\int_X (-\f) \o_{\f}^2.
$$
Since 
$\int_X (\f_j-\f) \o \wedge \o_{\f} \rightarrow 0$,
we infer
$$
\nu(X) \leq \limsup \int_X (-\f_j) \o_{\f_j}^2 \leq \int_X (-\f) \o_{\f}^2.
$$

Conversely if $\f \in {\mathcal E}^1(X,\o)$, the proof above shows that
$$
\limsup_{j \rightarrow +\infty} \int_X (-\f_j) \o_{\f_j}^2 \leq \int_X (-\f) \o_{\f}^2<+\infty,
$$
hence the sequence $(\int_X (-\f_j) \o_{\f_j}^2)$ is bounded
and $\f_j \o_{\f_j}^2 \rightarrow \f \o_{\f}^2$.
\end{proof}

We shall make constant use of theorem 3.1 in what follows. Indeed every function 
$\f \in {\mathcal E}^1$ is the decreasing limit of a sequence of smooth functions
$\f_j \in {\mathcal E}^1$ such that $\int (-\f_j) \o_{\f_j}^2$ is bounded. We can thus perform
integration by parts in the class ${\mathcal E}^1$ by working first with smooth approximants.
This will be implicit in our forthcoming computations.

\begin{pro}
\text{ }

1) Let  $\f \in {\mathcal E}^1(X,\o)$, $\f\leq 0$, then
$$
\int_X (-\f) \o^2 \leq \int_X (-\f) \o \wedge \o_{\f} \leq \int_X (-\f) \o_{\f}^2.
$$

2) If $\f,\p \in {\mathcal E}^1(X,\o)$ then 
${\mathcal E}^1(X,\o) \subset L^1(\o_{\f} \wedge \o_{\p})$.

3) The set ${\mathcal E}^1(X,\o)$ is a star-shaped convex.

4) Assume $\f,\p \in PSH(X,\o)$ with $\f \leq \p$. If 
$\f \in  {\mathcal E}^1(X,\o)$
then $\p \in {\mathcal E}^1(X,\o)$. In particular ${\mathcal E}^1(X,\o)$ is stable
under taking maximum.

5) Let $\f \in {\mathcal E}^1(X,\o)$, $\f \leq 0$. If $\f \in L^2(\o_{\f} \wedge \o)$, then
$d\f \wedge d^c \f \wedge \o_{\f}$ is a well defined positive measure whose total mass is
bounded by
$$
\int_X d\f \wedge d^c \f \wedge \o_{\f} \leq \int_X (-\f) \o_{\f}^2.
$$
\end{pro}

\begin{proof}
Let $\f \in {\mathcal E}^1(X,\o)$ with $\f \leq 0$. Assume first that $\f$ is smoooth. 
Observing that $d\f \wedge d^c \f \wedge \o_{\f}$ is a positive measure, we infer
from Stokes theorem that
$$
\int_X (-\f) \o \wedge \o_{\f} \leq \int_X (-\f) \o \wedge \o_{\f}+
\int_X d\f \wedge d^c \f \wedge \o_{\f}=\int_X (-\f) \o_{\f}^2.
$$
A similar use of Stokes theorem yields the second inequality. The general case now follows by
regularizing $\f$ and by using theorems 2.1 and 3.1. This proves 1).
\vskip.1cm

Let $\f,\p,u \in {\mathcal E}^1(X,\o)$. That $\o_{\f} \wedge \o_{\p}$ is a well defined 
probability measure follows from proposition 1.3. We can assume, w.l.o.g. that $\f,\p,u \leq 0$.
We are going to show that $u \in L^1(\o_{\f} \wedge \o_{\p})$ by proving
$$
\int_X (-u) \o_{\f} \wedge \o_{\p} \leq 6 M,
\; \; M:=\max\left( \int_X (-\f) \o_{\f}^2 \;
\int_X (-\p) \o_{\p}^2; \int_X (-u) \o_{u}^2 \right).
$$
By theorem 3.1 we can assume $\f,\p,u$ are smooth.
Observe first that by Stokes theorem,
\begin{eqnarray*}
&&\int (-u) \o_{\f} \wedge \o_{\p} = \int (-u) \o \wedge \o_{\p}+\int du \wedge d^c \f \wedge \o_{\p} \\
&\leq& \int (-u) \o \wedge \o_{\p}+\left( \int du \wedge d^c u \wedge \o_{\p} \right)^{1/2}
\cdot \left( \int d \f \wedge d^c \f \wedge \o_{\p} \right)^{1/2},
\end{eqnarray*}
where the last inequality follows from Cauchy-Schwarz inequality.
Now
$$
\int du \wedge d^c u \wedge \o_{\p}=\int (-u) dd^c u \wedge \o_{\p} \leq 
\int (-u) \o_u \wedge \o_{\p}.
$$
On the other hand it follows from 
Cauchy-Schwarz inequality again together with 3.2.1 that
$$
\int (-u) \o \wedge \o_{\p}=\int (-u) \o^2+\int du \wedge d^c \p \wedge \o \leq 2M.
$$
It suffices therefore to prove $\int_X (-\f) \o_{\f} \wedge \o_{\p} \leq 4M$.

We decompose again $\o_{\p}=\o+dd^c \p$ and integrate by parts to obtain
$$
\int (-\f) \o_{\f} \wedge \o_{\p} \leq M+M^{1/2} 
\left( \int (-\p) \o_{\p} \wedge \o_{\f} \right)^{1/2}.
$$
If $\int (-\p) \o_{\p} \wedge \o_{\f} \leq M$ we are done, otherwise this yields
\begin{eqnarray}
\int (-\f) \o_{\f} \wedge \o_{\p} \leq 2M^{1/2} 
\left(\int (-\p) \o_{\p} \wedge \o_{\f}\right)^{1/2}.
\end{eqnarray}
Similarly we obtain 
$$
\int (-\p) \o_{\p} \wedge \o_{\f} \leq M+M^{1/2} 
\left( \int (-\f) \o_{\f} \wedge \o_{\p} \right)^{1/2}.
$$
Either $\int (-\f) \o_{\f} \wedge \o_{\p} \leq M$ and we are done, or this yields
\begin{eqnarray}
\int (-\p) \o_{\p} \wedge \o_{\f} \leq 2M^{1/2} 
\left(\int (-\f) \o_{\f} \wedge \o_{\p}\right)^{1/2}.
\end{eqnarray}
Finally (4) and (5) yield the upper-bound $\int (-\f) \o_{\f} \wedge \o_{\p} \leq 4M$. 
This ends the proof of 2).
\vskip.1cm

Now 3) follows straightforwardly from 2), as in the proof of proposition 1.2.2. We turn to 4).
Assume $\f \leq \p \leq 0$ with $\f \in {\mathcal E}^1$. Then
\begin{eqnarray*}
0 \leq \int (-\p) \o_{\p}^2 &\leq& \int (-\f) \o_{\p}^2
 =\int (-\f) \o \wedge \o_{\p}+\int (-\p) dd^c \f \wedge \o_{\p} \\
&\leq&  \int (-\f) \o \wedge \o_{\p} +\int (-\f) \o_{\f} \wedge \o_{\p}.
\end{eqnarray*}
Going on integrating by parts, using $dd^c \f \leq \o_{\f}$, $-\p \leq -\f$ and 4.2.1, we end
up with
$$
0 \leq \int_X (-\p) \o_{\p}^2 \leq 4 \int_X (-\f) \o_{\f}^2,
$$
which proves 4).
\vskip.1cm

It remains to prove 5). Let $\f_j$ be smooth approximants of $\f$ and compute
$$
d\f_j \wedge d^c \f_j \wedge \o_{\f_j}=(-\f_j )\o_{\f_j}^2-(-\f_j) \o \wedge \o_{\f_j}
-dd^c\left( \frac{1}{2}\f_j^2 \o_{\f_j} \right)
$$
Thanks to theorem 3.1, all the terms on the right hand side converge if 
$\f \in {\mathcal E}^1$ and $\f \in L^2(\o \wedge \o_{\f})$. This shows that
the measure $d\f \wedge d^c \f \wedge \o_{\f}$ is well defined in this case.
Note that it is positive as a limit of positive measures. Moreover if $\f \leq 0$, then
$(-\f) \o \wedge \o_{\f} \geq 0$, hence
$$
\int_X d\f \wedge d^c \f \wedge \o_{\f} \leq \int_X (-\f) \o_{\f}^2.
$$
\end{proof}

Similar arguments as above now yield the following continuity result, 
whose proof is left to the reader.

\begin{thm}
The operator 
$$
(\f,\p,u) \mapsto u \o_{\f} \wedge \o_{\p}
$$
is continuous under decreasing sequences in ${\mathcal E}^1(X,\o)$.
\end{thm}

There is uniqueness of solutions to the Monge-Amp\`ere equation in the class ${\mathcal E}^1$
as the following result shows.

\begin{thm} [Uniqueness]
Let $\f,\p \in {\mathcal E}(X,\o)$ be such that $\o_{\f}^2 \equiv \o_{\p}^2$.
If $\f \in {\mathcal E}^1(X,\o)$, then
$\f-\p$ is constant.
\end{thm}

\begin{proof}
We first assume that $\p$ belongs to 
${\mathcal E}^1(X,\o)$ as well.
Set $f=\f-\p$ and $h=(\f+\p)/2 \in {\mathcal E}^1(X,\o)$ (proposition 3.2).
We assume w.l.o.g. $\f,\p \leq 0$.
We are going to prove that $\nabla f=0$ by establishing the following
inequality
\begin{eqnarray*}
\! \! \! \! \! \! \! \! \! \! \! \! \! \! \! \! 
\lefteqn{ {\bf (\dag)} \hskip.5cm
\int_X df \wedge d^c f \wedge \o \leq \int_X df \wedge d^c f \wedge \o_h} \\
&&\hskip2.5cm
+4 \left( \int_X df \wedge d^c f \wedge \o_h \right)^{1/2}
\cdot
\left( \int_X dh \wedge d^c h \wedge \o_h \right)^{1/2}.
\end{eqnarray*}
Observe that each integral on the right hand side is finite thanks to
proposition 3.2 if $\f \in L^2(\o_{\f} \wedge \o)$
and $\p \in L^2(\o_{\p} \wedge \o)$. 
Observe also that this yields the desired result in this case
since
$$
\int_X df \wedge d^c f \wedge \o_h=\int_X (-f) dd^c f \wedge \o_h
=\frac{1}{2}\int_X (-f) (\o_{\f}^2-\o_{\p}^2)=0,
$$
if $\o_{\f}^2 \equiv \o_{\p}^2$.When $\f \notin  L^2(\o_{\f} \wedge \o)$ or
$\p \notin  L^2(\o_{\p} \wedge \o)$,we use smooth approximants $\f_j,\p_j$ 
and observe that
$$
\int_X dh_j \wedge d^c h_j \wedge \o_{h_j} \leq \int_X (-h_j)\o_{h_j}^2 \leq 4 \int_X (-h) \o_h^2
$$
and 
$$
\int_X df_j \wedge d^c f_j \wedge \o_{h_j} =
\frac{1}{2}\int_X (\p_j-\f_j) \left(\o_{\f_j}^2-\o_{\p_j}^2 \right) \longrightarrow 0
$$
as follows from theorem 3.3.

We now establish $(\dag)$. Note that 
$df \wedge d^c f \wedge \o=df \wedge d^c f \wedge \o_h -df \wedge d^c f \wedge
dd^c h$, hence integrating by parts in the last term yields
$$
\int_X df \wedge d^c f \wedge \o=\int_X df \wedge d^c f \wedge \o_h
+\int_X df \wedge d^c h \wedge (\o_{\f}-\o_{\p}).
$$
Now it follows from Cauchy-Schwarz inequality that
$$
\left| \int_X df \wedge d^c h \wedge \o_{\f} \right|
\leq 2 \left(\int_X df \wedge d^c f \wedge \o_h \right)^{1/2}
\cdot
\left(\int_X dh \wedge d^c h \wedge \o_h \right)^{1/2}.
$$
A similar control on $\int_X df \wedge d^c f \wedge \o_{\p}$ yields $(\dag)$.

It remains to prove that $\p$ indeed belongs to 
${\mathcal E}^1(X,\o)$. 
We can assume $\p,\f \leq 0$.
Observe first that $\p \in L^1(\o_{\p} \wedge \o_{\f})$:
$$
0 \leq \int_X (-\p) \o_{\p} \wedge \o_{\f}=\int_X (\f-\p) \o_{\p} \wedge \o
+\int_X (-\f) \o_{\p}^2<+\infty
$$
since $\o_{\p}^2 \equiv \o_{\f}^2$ and $\f \in {\mathcal E}^1(X,\o)$. We infer
\begin{eqnarray*}
\lefteqn{\int_X (-\p) \o_{\p}^2 =
\int_X (-\p) \o_{\f}^2=\int_X (-\p) \o \wedge \o_{\f}
+\int_X (-\p) dd^c \f \wedge \o_{\f}} \\
&\leq& \int_X (-\p) \o \wedge \o_{\f}+
\left(\int_X d\p \wedge d^c \p \wedge \o_{\f} \right)^{1/2}
\cdot 
\left(\int_X d\f \wedge d^c \f \wedge \o_{\f} \right)^{1/2} \\
&\leq& \int_X (-\p) \o \wedge \o_{\f}+
\left(\int_X (-\p) \o_{\p} \wedge \o_{\f} \right)^{1/2}
\cdot 
\left(\int_X (-\f) \o_{\f}^2 \right)^{1/2}<+\infty.
\end{eqnarray*}

\end{proof}

\begin{rqe}
The idea of the proof of this uniqueness result goes back to E.Calabi [8]
who used the positivity of $df \wedge d^c f \wedge (\o_{\f}+\o_{\p})$ when
$\o_{\f},\o_{\p}$ are K\"ahler forms. The proof given above is a variation on
an argument of Z.Blocki [6] who proved the uniqueness in case $\f,\p$ are
bounded. 
\end{rqe}

The next lemma will be quite useful in section 6.
It gives a necessary condition for a probability measure to be the Monge-Amp\`ere 
of a function that belongs to the class ${\mathcal E}^1$.

\begin{lem}
Let $\mu$ be a probability measure on $X$. Then ${\mathcal E}^1(X,\o) \subset L^1(\mu)$
if and only if there exists $C_{\mu}>0$ such that for all functions 
$\f \in PSH(X,\o) \cap L^{\infty}(X)$ normalized by $\sup_X \f=-1$, one has
$$
0 \leq \int_X (-\f) d\mu \leq C_{\mu} \left( \int_X (-\f) \o_{\f}^2 \right)^{1/2}.
$$
\end{lem}

\begin{proof}
One direction is obvious. If there is such an inequality for all bounded $\o$-psh functions,
then the inclusion ${\mathcal E}^1 \subset L^1(\mu)$ follows from theorem 3.1.

Conversely assume the inequality is not satisfied, i.e. for all $j \in \N$, one can find
a bounded $\o$-psh function $\f_j$ such that $\sup_X \f_j=-1$ and 
$$
\int (-\f_j) d\mu > 2^j \left( \int (-\f_j) \o_{\f_j}^2 \right)^{1/2}. 
$$
Assume first
that the sequence $\int (-\f_j) \o_{\f_j}^2 \leq M$ is uniformly bounded from above.
We set then $\p:= \sum_{j \geq 1} 2^{-j} \f_j$. This is a 
well defined $\o$-psh function (as a decreasing
sequence of $\o$-psh functions which does not converge uniformly towards $-\infty$ thanks to the
normalisation $\sup_X \f_j=-1$). Now it follows from the estimate in the proof of
proposition 3.2.1 that 
$$
\int_X (-\p) \o_{\p}^2 \leq 6 \sup_{j \geq 1} \int_X (-\f_j) \o_{\f_j}^2 <+\infty
$$
hence $\p \in {\mathcal E}^1$, while $\int (-\p) d\mu=+\infty$ by the monotone convergence theorem.

Assume now $M_j:=\int (-\f_j) \o_{\f_j}^2 \rightarrow +\infty$. We set $\p_j:=\e_j \f_j$, where
$\e_j=M_j^{-1/2}$ is chosen so that $\int (-\p_j) \o_{\p_j}^2 \leq M$ is uniformly bounded. Indeed
a straightforward computation yields
\begin{eqnarray*}
\int (-\p_j) \o_{\p_j}^2 &\leq& \e_j \int (-\f_j)\o^2 +2\e_j^2 \int (-\f_j)\o \wedge \o_{\f_j}
+\e_j^3 \int (-\f_j) \o_{\f_j}^2 \\
&\preceq& \e_j^2 \int (-\f_j) \o_{\f_j}^2 \leq M,
\end{eqnarray*}
because $\int (-\f_j) \o^2$ is uniformly bounded since $\sup_X \f_j=-1$ and 
\linebreak[4]
$\int (-\f_j) \o \wedge \o_{\f_j} \leq \int (-\f_j) \o_{\f_j}^2$.
We set now $\p:=\sum 2^{-j} \p_j$. This is a well defined function in ${\mathcal E}^1$
such that
$$
\int (-\p_j) d\mu=\e_j \int(-\f_j) d\mu >2^j \e_j \left( \int (-\f_j) \o_{\f_j}^2 \right)^{1/2}=2^j.
$$
Thus $\int (-\p) d\mu=+\infty$, so that ${\mathcal E}^1$ is not included in $L^1(\mu)$.
\end{proof}

When $\mu=\o_{\p}^2$ is the Monge-Amp\`ere of a function $\p \in {\mathcal E}^1(X,\o)$,
it follows from proposition 3.2 that ${\mathcal E}^1(X,\o) \subset L^1(\mu)$,
hence there exists $C_{\p}>0$ such that for all functions
$\f \in PSH(X,\o) \cap L^{\infty}(X)$ normalized by $\sup_X \f=-1$, one has
$$
0 \leq \int_X (-\f) \o_{\p}^2 \leq C_{\p} \left[ \int_X (-\f) \o_{\f}^2 \right]^{1/2}.
$$

\section{The classes ${\mathcal E}^p(X,\o)$ }

In this section we fix a real number $p \geq 1$.

\begin{defi}
We let ${\mathcal E}^p(X,\o)$ denote the set of functions $\f \in {\mathcal E}(X,\o)$
such that there exists a sequence $\f_j \in PSH(X,\o) \cap L^{\infty}(X)$ with
$$
\f_j \searrow \f \; \text{ and } \; 
\sup_{j \in \N} \left( \int_X |\f_j|^p \o_{\f_j}^2\right) < +\infty.
$$
\end{defi}

When no confusion can arise, we shall simply denote 
${\mathcal E}^p(X,\o)$ by ${\mathcal E}^p$. Similar classes were introduced by
U.Cegrell in the local context [9] as generalizations of the classical
notion of subharmonic functions of finite energy. Observe that
${\mathcal E}(X,\o)={\mathcal E}^0(X,\o) \supset {\mathcal E}^p(X,\o) \supset
{\mathcal E}^q(X,\o)$ whenever $p \leq q$.
When $p=1$ this definition is equivalent to the one we gave in section 3, thanks
to theorem 3.1. When $p>1$ we of course get the inclusion
$$
{\mathcal E}^p(X,\o) \subset \left\{ \f \in {\mathcal E}(X,\o) \, / \,
\f \in L^p(\o_{\f}^2) \right \},
$$
however the reverse inclusion is not clear: we don't know how to produce a decreasing
sequence with uniformly bounded energies.
Indeed a delicate point in the analysis of the classes $
{\mathcal E}^p$ is that we don't know
if a continuity result similar to theorem 3.1 still holds. We shall prove a
weaker property in theorem 4.4 below, but we need first
to establish some useful inequalities.

\begin{lem}
Let $\f,\p \in PSH(X,\o) \cap L^{\infty}(X)$ 
with $\f \leq \p \leq 0$. Then

1) $0 \leq \int_X (-\f)^p \o^2 \leq \int_X (-\f)^p \o \wedge \o_{\f}
\leq \int_X (-\f)^p \o_{\f}^2$.

2) $0 \leq \int_X (-\p)^p \o \wedge \o_{\p} \leq (p+1) 
\int_X (-\f)^p \o \wedge \o_{\f}$.

3) $0 \leq \int_X (-\p)^p \o_{\p}^2 \leq (p+1)^2 \int_X (-\f)^p \o_{\f}^2.$
\end{lem}

\begin{proof}
It suffices to establish these inequalities when $\f,\p$
are smooth. Indeed one can approximate $\f$, $\p$ by decreasing sequences of
smooth $\o$-psh functions and then use classical continuity results of
E.Bedford and A.Taylor [5].

In the sequel we thus assume $\f,\p$ are smooth. Observe that 
$d \f \wedge d^c \f \wedge T$ is a positive measure whenever $T$ is a smooth $(1,1)$-form.
It follows therefore from Stokes theorem that
$$
\int_X (-\f)^p \o_{\f}^2 =\int_X (-\f)^p \o \wedge \o_{\f}
+p\int_X (-\f)^{p-1} d\f \wedge d^c\f \wedge \o_{\f}
\geq  \int_X (-\f)^p \o \wedge \o_{\f}
$$
and, similarly,
$$
\int_X (-\f)^p \o \wedge \o_{\f} =\int_X (-\f)^p \o^2
+p\int_X (-\f)^{p-1} d\f \wedge d^c\f \wedge \o
\geq  \int_X (-\f)^p \o^2.
$$
This proves 1).
Now observe that
$$
-dd^c (-\f)^p=p(-\f)^{p-1} dd^c \f -p(p-1)(-\f)^{p-2}d\f \wedge d^c \f 
\leq p(-\f)^{p-1} \o_{\f}.
$$
It follows therefore from Stokes theorem that
\begin{eqnarray*}
\int (-\p)^p \o \wedge \o_{\p} &\leq& \!\!\!\!\! \int (-\f)^p \o \wedge \o_{\p}
=\int (-\f)^p \o^2+\int (-\p) [-dd^c (-\f)^p] \wedge \o \\
&\leq& \! \! \! \! \! \int (-\f)^p \o^2 +p \int (-\f)^p \o \wedge \o_{\f}
\leq (p+1) \int (-\f)^p \o \wedge \o_{\f},
\end{eqnarray*}
which proves 2).
The proof of the third inequality is similar and is left to the reader. 
\end{proof}

\begin{cor}
Let $\f \in {\mathcal E}^p(X,\o)$. Let $\p_j \in PSH(X,\o) \cap L^{\infty}(X)$ 
be {\bf any} sequence decreasing towards $\f$. Then
$$
\sup_{j \in \N} \left( \int_X |\p_j|^p \o_{\p_j}^2 \right)<+\infty.
$$
\end{cor}

\begin{proof} 
Let $\f_j$ be a sequence of bounded $\o$-psh functions which decreases 
towards $\f$ and with bounded energies. 
We can assume w.l.o.g. that $\f_j,\p_k \leq 0$ for all $j,k$.
We fix $k$ and consider the sequence
$\Phi_j:=\max(\f_j,\p_k) \in PSH(X,\o) \cap L^{\infty}$.
This is a sequence of uniformly bounded $\o$-psh functions such that
$\Phi_j \downarrow \p_k$ as $j \rightarrow +\infty$.
It follows therefore from classical continuity results of Bedford and Taylor
that $(-\Phi_j)^p \o_{\Phi_j}^2\rightarrow (-\p_k)^p \o_{\p_k}^2$.
Now $\Phi_j \geq \f_j$, so it follows from lemma 4.2 that
$$
\int_X (-\p_k) \o_{\p_k}^2=\lim_{j \rightarrow +\infty} \int_X (-\Phi_j)^p \o_{\phi_j}^2
\leq (p+1)^2 \sup_{j \geq 1} \int_X (-\f_j)^p \o_{\f_j}^2.
$$
\end{proof}

Recall that a sequence of $\o$-psh functions $(\f_j)$ converges in capacity
towards $\f \in PSH(X,\o)$ if for all $\e>0$,
$$
Cap_{\o} \left( |\f_j-\f|>\e \right) \longrightarrow 0.
$$
Following [24] we now show that convergence in capacity implies convergence
of Monge-Amp\`ere operators. 

\begin{thm}
Let $\f \in {\mathcal E}^p(X,\o)$ and let $(\f_j)$ be a sequence of $\o$-psh functions
that converges in capacity towards $\f$.
If $\sup (\int |\f_j|^p \o_{\f_j}^2)<+\infty$, then
for all $q<p$, $\f_j \in {\mathcal E}^q(X,\o)$ and
$$
(-\f_j)^q \o_{\f_j}^2 \longrightarrow (-\f)^q \o_{\f}^2.
$$
\end{thm}

\begin{proof}
The result is true, even with $q=p$, when the functions $\f_j$ are uniformly
bounded (see [24]). We are going to reduce to that case by
considering 
$$
\f^{(k)}:=\max(\f,-k) \, \text{ and } \; 
\f_j^{(k)}:=\max(\f_j,-k).
$$
Let $\chi$ be a test function. We want to show that
$<(-\f_j)^q\o_{\f_j}^2,\chi> \rightarrow <(-\f)^q\o_{\f}^2,\chi>$.
Decomposing in the obvious way, this boils down to establish 
good upper bounds on 
$$
\int_{(\f_j \leq -k)} (-\f_j)^q \o_{\f_j}^2, \; \; 
\int_{(\f_j \leq -k)} (-\f_j^{(k)})^q \o_{\f_j^{(k)}}^2, \; \; 
\int_{(\f \leq -k)} (-\f^{(k)})^q \o_{\f^{(k)}}^2
$$
that are uniform in $j$. By Chebyshev inequality we get
$$
\int_{(\f_j \leq -k)} (-\f_j)^q \o_{\f_j}^2
\leq \frac{1}{k^{p-q}} \int_X (-\f_j)^q \o_{\f_j}^2
\leq \frac{C}{k^{p-q}}.
$$
We make a similar use of Chebyshev inequality on the two other integrals.
Now it follows from lemma 4.2 and corollary 4.3 that
$$
\int_X \left( -\f_j^{(k)} \right)^q \o_{\f_j^{(k)}}^2, \; 
\int_X \left( -\f^{(k)} \right)^q \o_{\f^{(k)}}^2 
\leq (p+1)^2 \int_X (-\f_j)^q \o_{\f_j}^2.
$$
\end{proof}

Note that if $\f_j \in PSH(X,\o) \cap L^{\infty}(X)$ decreases towards $\f$,
then $\f_j$ converges towards $\f$ in capacity with
$\sup \int |\f_j|^p \o_{\f_j}^2<+\infty$ (by corollary 4.3). It follows
therefore from theorem 4.4 that lemma 4.2 holds whenever
$\f,\p$ are in ${\mathcal E}^p(X,\o)$.

\begin{cor}
Lemma 4.2 holds with $\f , \p \in {\mathcal E}^p(X,\o)$.
In particular ${\mathcal E}^p(X,\o)$ is stable under taking maximum.
\end{cor}

\begin{proof}
Theorem 4.4 and lemma 4.2 show that lemma 4.2 holds for any exponent
$q<p$ and for any fixed function $\f \in {\mathcal E}^p(X,\o)$.
Now letting $q$ increase towards $p$ yields the conclusion thanks
to the monotone convergence theorem.
This shows in particular that if $\f \in {\mathcal E}^p(X,\o))$
and $\f \leq \p$, $\p \in PSH(X,\o)$, then $\p \in {\mathcal E}^p(X,\o))$.
Thus ${\mathcal E}^p(X,\o))$ is stable under taking maximum.
\end{proof}

\begin{pro}
\text{ }

1) If $\f,\p \in {\mathcal E}^p(X,\o)$, then the probability measure $\o_{\f} \wedge \o_{\p}$
satisfies
$$
{\mathcal E}^p(X,\o) \subset L^p(\o_{\f} \wedge \o_{\p}).
$$

2) The set ${\mathcal E}^p(X,\o)$ is a star-shaped convex.

3) If $\f \in {\mathcal E}^p(X,\o)$ is such that
$\f \in L^{p+1}(\o \wedge \o_{\f})$,  then the positive measure 
$(-\f)^{p-1}d\f \wedge d^c \f \wedge \o_{\f}$ is well defined and has total
mass
$$
\int_X (-\f)^{p-1}d\f \wedge d^c \f \wedge \o_{\f} \leq \frac{1}{p}
\int_X (-\f)^p \o_{\f}^2.
$$
\end{pro}

\begin{proof}
Assume $\f,\p \in {\mathcal E}^p$.
In particular $\f,\p \in {\mathcal E}$ so $\o_{\f} \wedge \o_{\p}$
is a well defined probability measure. Let $u \in {\mathcal E}^p$; we assume
w.l.o.g. $u,\f,\p \leq 0$. We are going to give an upper bound on
$\int (-u)^p \o_{\f} \wedge \o_{\p}$ which only depends on
$\int (-u)^p \o_{u}^2$, $\int (-\f)^p \o_{\f}^2$, $\int (-\p)^p \o_{\p}^2$.
Using theorem 4.4 it suffices to establish this bound when $u,\f,\p$ are smooth.
Recall that $-dd^c(-u)^p \leq p(-u)^{p-1} \o_u$ and
$\int (-u)^p \o \wedge \o_{\p} \leq \int (-u)^p \o_u \wedge \o_{\p}$,
so it follows from Stokes theorem and H\"older inequality that
\begin{eqnarray*}
\int && \!\!\!\!\!\!\!\!\!\!\!\!\!\!\!\!\!
(-u)^p \o_{\f} \wedge \o_{\p} = \int (-u)^p \o \wedge \o_{\p}
+\int (-\f) [-dd^c (-u)^p] \wedge \o_{\p} \\
&\leq& \!\!\!\!\!\!(p+1) \left[ \int (-u)^p \o_u \wedge \o_{\p} \right]^{1-1/p}
\left[ \max \left\{ \int (-u)^p \o_u \wedge \o_{\p}; 
\int (-\f)^p \o_u \wedge \o_{\p} \right\} \right]^{1/p}.
\end{eqnarray*}
Reversing the roles of $u,\f$ we infer
$$
\int (-u)^p \o_{\f} \wedge \o_{\p} \leq (p+1)^{\frac{p}{p-1}}
\max \left[ \int (-u)^p \o_u \wedge \o_{\p}; 
\int (-\f)^p \o_{\f} \wedge \o_{\p}  \right].
$$
Thus it suffices to get an upper bound when $u=\f$.
We set 
$$
M:=\max \left[ \int (-\f)^p \o_{\f}^2 ; \int (-\p)^p \o_{\p}^2 \right].
$$
We use our last inequality with $u=\p$ to obtain
$$
\int (-\p)^p \o_{\f} \wedge \o_{\p} \leq (p+1)M^{1-1/p} 
\left[ \max \left\{ M; \int (-\f)^p \o_{\f} \wedge \o_{\p} \right\} \right]^{1/p}.
$$
Reversing the roles of $\f$ and $\p$, this yields 
$$
\int (-\f)^p \o_{\f} \wedge \o_{\p} \leq (p+1)^{\frac{p}{p-1}} M,
$$
which proves 1).

As in the proof of proposition 1.3.2,we infer straightforwardly from 1) that 
${\mathcal E}^p$ is a star-shaped convex. 
When $\f \in L^{p+1}(\o \wedge \o_{\f})$, the positive current
$(-\f)^{p+1} \o_{\f}$ is well defined and we compute
$$
\frac{1}{(p+1)}
dd^c [(-\f)^{p+1} \o_{\f}]=p (-\f)^{p-1} d\f \wedge d^c \f \wedge \o_{\f}
- (-\f)^p \o_{\f}^2+ (-\f)^p \o_{\f} \wedge \o.
$$
This allows us to define
$(-\f)^{p-1} d\f \wedge d^c \f \wedge \o_{\f}$ as soon as $p>0$.
It then follows from Stokes theorem that
$$
\int_X (-\f)^p \o_{\f}^2=p \int_X (-\f)^{p-1} d\f \wedge d^c \f \wedge \o_{\f}
+\int_X (-\f)^p \o_{\f} \wedge \o,
$$
which yields the desired upper bound 3).
\end{proof}

The next lemma is analogous to lemma 3.6; it will be quite useful
in section 5 as it gives a necessary condition for a probability measure to be
the Monge-Amp\`ere of a function that belongs to ${\mathcal E}^p$.

\begin{lem}
Let $\mu$ be a probability measure on $X$. Then ${\mathcal E}^p(X,\o) \subset L^p(\mu)$
if and only if there exists $C>0$ such that for all functions $\f \in PSH(X,\o) \cap L^{\infty}(X)$
normalized by $\sup_X \f=-1$, one has
$$
0 \leq \int_X (-\f)^p d\mu \leq C \left( \int_X (-\f)^p \o_{\f}^2 \right)^{\frac{p}{p+1}}.
$$
\end{lem}

The proof is very similar to that of lemmas 3.6. We leave it to the reader.

We end this section with a result that will be crucial when solving Monge-Amp\`ere
equations in the next section.

\begin{thm}
Let $\f_j \leq 0$ be a sequence of functions in ${\mathcal E}^p(X,\o)$ such that
$$
\sup_{j \in \N} \int_X (-\f_j)^p \o_{\f_j}^2 <+\infty.
$$
Assume $\f_j \rightarrow \f$ in $L^1(X)$. Then $\f \in {\mathcal E}^p(X,\o)$.

If moreover $\int_X |\f_j-\f| \o_{\f_j}^2 \rightarrow 0$, then
$\o_{\f_j}^2 \rightarrow \o_{\f}^2$.
\end{thm}

\begin{proof}
The bounded $\o$-psh functions 
$$
\p_j:=\left( \sup_{k \geq j} \max [\f_k,-k] \right)^*
$$
decrease towards $\f$. Since $\p_j \geq \f_j$, it follows from corollary 4.5 that
$\sup_j \int (-\p_j)^p \o_{\p_j}^2<+\infty$, therefore $\f \in {\mathcal E}^p(X,\o)$.

We assume now that $\int |\f_j-\f| \o_{\f_j}^2 \rightarrow 0$. Passing to a
subsequence if necessary, we can assume $\int |\f_j-\f| \o_{\f_j}^2 \leq 1/j^2$.
Consider
$$
\Phi_j:=\max (\f_j,\f-1/j) \in {\mathcal E}^p(X,\o).
$$
It follows from Hartogs' lemma that $\Phi_j \rightarrow \f$ in capacity, hence
$\o_{\Phi_j}^2 \rightarrow \o_{\f}^2$ by theorem 4.4.
Thus we need to compare $\o_{\Phi_j}^2$ and $\o_{\f_j}^2$. It follows from a classical
inequality of J.-P. Demailly [11] that
$$
\o_{\Phi_j}^2 \geq \1_{\{ \f_j \geq \f-1/j \}} \cdot \o_{\f_j}^2.
$$
Let $E_j$ denote the set $X \setminus \{\f_j \geq \f-1/j\}$, i.e.
$E_j=\{\f-\f_j>1/j\}$. Our assumption implies that 
$\1_{E_j} \o_{\f_j}^2 \rightarrow 0$, indeed
$$
0 \leq \int_{E_j} \o_{\f_j}^2 \leq j \int_X |\f-\f_j| \o_{\f_j}^2 \leq \frac{1}{j}.
$$
Therefore $0 \leq \o_{\f_j}^2 \leq \o_{\Phi_j}^2+o(1)$, hence
$\o_{\f}^2=\lim \o_{\f_j}^2$.
\end{proof}

\section{Range of the complex Monge-Amp\`ere operator}

In this section we prove our main result. This is the following

\begin{thm}
Let $\mu$ be a probability measure on $X$ and $p \geq 1$.

Then there exists a unique function $\p \in {\mathcal E}^p(X,\o)$ such that
$$
\mu=\o_{\p}^2 \; \; \text{ and } \; \; \sup_X \p=-1
$$
if and only if ${\mathcal E}^p(X,\o) \subset L^p(\mu)$.
\end{thm}

This result follows straightforwardly from lemmas 3.6, 4.7 together with the 
following theorem (applied with $\a=p/(p+1)$).

\begin{thm} 
Fix $p \geq 1$, $0 <\a<1$ and $C>0$.
Let $\mu$ be a probability measure such that for all functions
$\f \in PSH(X,\o) \cap L^{\infty}(X)$ with $\sup_X \f=-1$, one has
$$
0 \leq \int_X (-\f)^p d\mu \leq C \left( \int_X (-\f)^p \o_{\f}^2 \right)^{\a}.
$$
Then there exists a unique function $\p \in {\mathcal E}^p(X,\o)$ s.t.
$\mu=\o_{\p}^2$, $\sup_X \p=-1$.
\end{thm}

The uniqueness of the solution $\p$, once normalized by $\sup_X \p=-1$, follows from
theorem 3.4. The proof of the existence will occupy the rest of this section.
The strategy of the proof is as follows:
\begin{itemize}

\item We approximate $\mu$ by smooth probability volume forms $\mu_j$
 using local convolutions and a partition of unity.

\item We invoke Aubin-Yau's solution of the Calabi conjecture to find
 uniquely determined $\o$-psh functions $\f_j$ such that
 $\mu_j=\o_{\f_j}^2$ and $\sup_X \f_j=-1$.

\item Since $\o$-psh functions $\f$ normalized by $\sup_X \f=-1$ form a compact subset
 of $L^1(X)$, we can assume that $\f_j \rightarrow \f$ in $L^1(X)$.

\item The integrability condition on $\mu$ guarantees $\sup_j \int (-\f_j)^p \o_{\f_j}^2 <+\infty$,
 hence yields $\f \in {\mathcal E}^p(X,\o)$.

\item The delicate point is then to show that $\o_{\f_j}^2 \rightarrow \o_{\f}^2$. This is done
 by showing that $\int |\f_j-\f| d\mu_j \rightarrow 0$ and invoking theorem 4.8.
 here we use the integrability assumption again with $p>1$ in order to show first that
 $\int \f_j d\mu \rightarrow \int \f d\mu$.

\item The case $p=1$ deserves special treatment. 
\end{itemize}

Here follow the technical details. Let $\{U_i\}$ be a finite covering of $X$ by open sets 
$U_i$ which are biholomorphic to the unit ball of $\C^2$. In each $U_i$ we let
$\mu_{\e}^{U_i}:=\mu_{|U_i} * \chi_{\e}$ denote local regularization of $\mu_{|U_i}$ by means
of convolution with radial nonnegative smooth approximations $\chi_{\e}$ of the Dirac mass.
Let $\{\theta_i\}$ be a partition of unity subordinate to $\{U_i\}$ and set
$$
\mu_j:=c_j \left[ \sum_i \theta_i \mu_{\e_j}^{U_i}+\e_j \o^n \right],
$$
where $\e_j \searrow 0$ and $c_j \nearrow 1$ is chosen so that $\mu_j(X)=1$.
Thus the $\mu_j$'s are smooth probability volume forms which converge weakly towards $\mu$.
It follows from the solution of the Calabi conjecture [1], [25], that there exists a unique
function $\f_j \in PSH(X,\o) \cap {\mathcal C}^{\infty}(X)$ such that
$$
\mu_j=\o_{\f_j}^2 \; \; \text{ and } \; \; \sup_X \f_j=-1.
$$
Recall from proposition 1.7 in [16] that ${\mathcal F}:=\{ \f \in PSH(X,\o) \, / \, \sup_X \f=-1\}$
is a compact subset of $L^1(X)$. Passing to a subsequence if necessary, we can therefore assume
$\f_j \rightarrow \f$ in $L^1(X)$, where $\f \in PSH(X,\o)$ with $\sup_X \f=-1$.

\begin{lem}
There exists $C>1$ such that for all $j \in \N$,
$$
\int_X (-\f_j)^p \o_{\f_j}^2 \leq C \int_X (-\f_j)^p d\mu \leq C^2.
$$
In particular $\f \in {\mathcal E}^p(X,\o)$.
\end{lem}

\begin{proof}
Since $c_j \rightarrow 1$ and $\e_j \rightarrow 0$, we can write
$$
\int_X (-\f_j)^p \o_{\f_j}^2=\sum_i \int_X \theta_i (-\f_j)^p d\mu_{\e_j}^{U_i}+o(1),
$$
where 
$$
\int_X \theta_i (-\f_j)^p d \mu_{\e_j}^{U_i} =\int_{U_i} \theta_i * \chi_{\e_j} (-\f_j*\chi_{\e_j})^p d\mu
+o(1).
$$
Now $\f_j=u_j^i-\g_i$ in $U_i$, where $\g_i$ is a smooth local potential of $\o$ in $U_i$ and
$u_j^i$ is psh in $U_i$. Therefore
$-u_i^i*\chi_{\e_j} \leq -u_j^i$, while $\g_i * \chi_{\e_j}$ and $\theta_i * \chi_{\e_j}$ both converge 
uniformly towards $\g_i$ and $\theta_i$. We infer
$$
\int_{U_i} \theta_i * \chi_{\e_j} (-\f_j*\chi_{\e_j})^p d\mu \leq 
\int_{U_i} \theta_i (-\f_j )^p d\mu +o(1)
$$
hence 
$$
\int (-\f_j)^p d\mu_j= \int (-\f_j)^p \o_{\f_j}^2 \leq \int (-\f_j)^pd\mu+o(1).
$$
It follows now from the integrability assumption we made on $\mu$ that
$$
\sup_{j \in \N} \left( \int_X (-\f_j)^p \o_{\f_j}^2 \right)^{1-\a}<+\infty,
$$
which yields the lemma.
\end{proof}

We now would like to apply theorem 4.8 to insure that $\mu=\o_{\f}^2$. For this we need to assume $p>1$
in order to use the following:

\begin{lem}
Assume $p>1$. Then
$$
\int_X \f_j d\mu \rightarrow \int_X \f d\mu \; \; 
\text{ and } \; \; 
\int_X |\f_j-\f|d\mu_j \rightarrow 0.
$$
\end{lem}

\begin{proof}
When the $\f_j$'s are uniformly bounded, the first convergence follows from standard arguments (see
[9]). Set
$$
\f_j^{(k)}:=\max(\f_j,-k) \; \; \text{ and } \; \; \f^{(k)}:=\max(\f,-k).
$$
We will be done with the first convergence if we can show that $\int |\f_j^{(k)}-\f_j|d\mu \rightarrow 0$
uniformly in $j$ as $k \rightarrow +\infty$. This is where we use the assumption $p>1$. Namely
$$
0 \leq \int_X  |\f_j^{(k)}-\f_j|d\mu \leq 
2 \int_{(\f_j<-k)} (-\f_j)d\mu \leq \frac{2}{k^{p-1}} \int_X (-\f_j)^p d\mu \leq \frac{C}{k^{p-1}}.
$$
It remains to prove a similar convergence when $\mu $ is replaced by $\mu_j$.
It actually suffices to consider the case of measures
$\mu_j^U:=\mu_{|U}*\chi_{\e_j}$. Now
$$
\int_U |\f_j-\f|d\mu_j^U=\int_U \left(\int_U |u_j(\zeta)-u(\zeta)| 
\chi_{\e_j}(z-\zeta) d\lambda(\zeta) \right) d\mu(z),
$$
where as above, $u_j$, $u$ are psh functions in $U$ such that $\f_j=u_j-\g$ and  $\f=u-\g$ in $U$,
$\g$ is a local potential of $\o$ in $U$ and $d\lambda$ denotes the Lebesgue measure in $U$.
The lemma will be proved if we can show that $\int w_j d\mu \rightarrow 0$, where
$$
w_j(z):=\int_U |u_j(\zeta)-u(\zeta)| \chi_{\e_j}(z-\zeta) d\lambda(\zeta).
$$
Define $\tilde{u}_j:=(\sup_{k \geq j} u_k )^*$. This is a sequence of psh functions in $U$ which decrease
towards $u$. Observe that $\tilde{u}_j \geq \max(u,u_j)$ so that
$$
w_j \leq 2 \tilde{u}_j * \chi_{\e_j} -u * \chi_{\e_j}-u_j* \chi_{\e_j}
\leq 2(\tilde{u}_j*\chi_{\e_j}-u)+(\f-\f_j).
$$
It follows from the monotone convergence theorem that
$\int (\tilde{u}_j*\chi_{\e_j}-u) d\mu \rightarrow 0$, while 
$\int (\f_j-\f)d\mu \rightarrow 0$ by the first part of lemma.
Therefore $\int w_j d\mu \rightarrow 0$ and we are done.
\end{proof}

It follows from previous lemma and theorem 4.6 that $\mu=\o_{\f}^2$. This proves theorem 5.2 when $p>1$.
Assume now $p=1$. Following an idea of U.Cegrell [9] we consider the set
${\mathcal C}$ of probability measures $\nu$ such that for all $\f \in
PSH(X,\o) \cap L^{\infty} (X,\o)$,
$\sup_X \f=-1$, one has
$$
0 \leq \int_X (-\f)^2 d\nu \leq C_0 \left( \int_X (-\f)^2 \o_{\f}^2 \right)^{1/2},
$$
where $C_0=C_0(X,\o)$ is the constant given by lemma 5.5 below.
The set ${\mathcal C}$ is a convex compact set of probability measures which contains all measures
$\o_u^2$, where $u \in PSH(X,\o)$ is such that $-1 \leq u \leq 0$: this is the contents of lemma 5.5 below.
It follows from a generalization of Radon-Nikodym theorem [20] that one can decompose
$$
\mu=f \nu+\sigma, \; \; 
\text{ where } \nu \in {\mathcal C}, \; f \in L^1(\nu) \;  \text{ and } \sigma \in {\mathcal C}^{\perp}.
$$
Now $\sigma$ is carried by a pluripolar set because ${\mathcal C}$ contains all the measures
$\o_u^2$, $-1 \leq u \leq 0$, and $\mu$ does not charge pluripolar sets because
${\mathcal E}^1(X,\o) \subset L^1(\mu)$, thus $\sigma=0$. Consider
$$
\mu_j:= \d_j \min(f,j) \nu,
$$
where $\d_j \searrow 1$ so that $\mu_j(X)=1$. It follows from theorem 5.2 (case $p=2$) that
there exists a unique $\f_j \in {\mathcal E}^2(X,\o)$ with $\sup_X \f_j =-1$ and $\mu_j=\o_{\f_j}^2$.
We can assume $\f_j \rightarrow \f$ in $L^1(X)$.
Now
$$
\int (-\f_j) \o_{\f_j}^2 \leq \d_j \int (-\f_j) d\mu \leq C \d_j \left[ (-\f_j) \o_{\f_j}^2 \right]^{\a},
$$
so that $\sup_j \int (-\f_j) \o_{\f_j}^2<+\infty$ hence $\f \in {\mathcal E}^1(X,\o)$
(see theorem 4.8).
We set
$$
\Phi_j:= (\sup_{k \geq j} \f_k)^* \; \; \text{ and } F_j:=\inf_{k \geq j} \d_k \min(f,k).
$$
Clearly $\Phi_j \in {\mathcal E}^1(X,\o)$ with $\Phi_j \searrow \f$ and $F_j \nearrow f$. It follows from a classical
inequality of J.-P.Demailly [11] that
$$
\o_{\Phi_j}^2 \geq F_j \nu.
$$
We infer $\o_{\f}^2 \geq \mu$, whence equality since these are both probability measures.
This finishes the proof of theorem 5.2.
\vskip 0.1 cm

\begin{lem}
There exists $C_0=C_0(X,\o) >1$ such that for all $\f \in PSH(X,\o) \cap L^{\infty}(X)$, $\sup_X \f=-1$,
and for all $u \in PSH(X,\o) \cap L^{\infty}(X)$, $-1 \leq u \leq 0$, one has
$$
0 \leq \int_X (-\f)^2 \o_u^2 \leq C_0 \left( \int_X (-\f)^2 \o_{\f}^2 \right)^{1/2}.
$$
\end{lem}

\begin{proof}
It follows from Stokes theorem and Cauchy-Schwarz inequality that
\begin{eqnarray*}
\int (-\f)^2 \o_{u}^2 -\int (-\f)^2 \o^2 &\leq & 2 \int (-\f) (-u) \o_{\f} \wedge [\o+ \o_u] \\
&\leq & 2 \sqrt{2} \left( \int (-\f)^2 \o_{\f} \wedge [\o+\o_u] \right)^{1/2} \\
& \leq& 2 \sqrt{2} \left( 2 \int (-\f)^2 \o_{\f}^2+\int(-\f)^2 \o_{\f} \wedge dd^c u \right)^{1/2}.
\end{eqnarray*}
Now 
$\int (-\f)^2 \o_{\f} \wedge dd^c u \leq 
2 \int (-u)(-\f) \o_{\f}^2 \leq 2 \int (-\f)^2 \o_{\f}^2$ and
$\int (-\f)^2 \o^2$ is bounded from above by a uniform constant that only depends on $X,\o$, since
we have normalized $\f$ by $\sup_X \f=-1$. The lemma follows.
\end{proof}

\section{Capacity of sublevel sets}

In this section we want to connect the condition $\f \in {\mathcal E}^p(X,\o)$
to the size of the sublevel sets $(\f<-t)$ measured by the complex 
Monge-Amp\`ere capacity $Cap_{\o}$: the smaller $(\f<-t)$, the better the exponent
$p$ for which $\f \in {\mathcal E}^p(X,\o)$. This will allow us to give several examples
of unbounded functions $\f \in {\mathcal E}^p(X,\o)$, as well as examples
of probability measures $\mu$ such that ${\mathcal E}^p(X,\o) \subset L^p(\mu)$.
The basic tool to establish this connection is the comparison principle
which we now recall (see [19] for a proof).

\begin{pro}
Let $\f,\p \in PSH(X,\o) \cap L^{\infty}(X)$. Then
$$
\int_{(\f<\p)} \o_{\p}^2 \leq \int_{(\f < \p)} \o_{\f}^2.
$$
\end{pro}

As earlier $p$ denotes a real number $\geq 1$ and 
$\o$-psh functions $\f$ are normalized so that $\sup_X \f=-1$ (unless otherwise specified).
We set
$$
e_p(\f):=\int_X (-\f)^p \o_{\f}^2+2 \int_X (-\f)^{p+1} \o \wedge \o_{\f} 
+\int_X (-\f)^{p+2} \o^2.
$$
This is a well defined quantity as soon as $\f \in {\mathcal E}^p(X,\o)$ and
$\f \in L^{p+1}(\o \wedge \o_{\f})$, for instance when $\f \in {\mathcal E}^{p+1}(X,\o)$.
In the sequel we shall say that a function $h \geq 1$ belongs to $L^q(Cap_{\o})$ if the
following integral converges
$$
0 \leq \int_X h^q dCap_{\o}:=q \int_1^{+\infty} t^{q-1} Cap_{\o}(h>t) dt<+\infty.
$$
The following lemma establishes the basic connection between capacity of sublevel sets
and the energy $e_p(\f)$.

\begin{lem}
Assume $\f \in {\mathcal E}^p(X,\o) \cap L^{p+1}(\o \wedge \o_{\f})$ is normalized
so that $\sup_X \f=-1$. Then
$$
\frac{p+2}{p} \int_X (-\f)^p \o_{\f}^2 \leq
\int_X (-\f)^{p+2} dCap_{\o} \leq
2^{p+2} e_p(\f).
$$

In particular if $\p \in {\mathcal E}^p(X,\o)$, then there exists $C_{\p}>0$ s.t.
for all $t \geq 1$,
$
Cap_{\o}(\p <-t) \leq C_{\p}/t^{p+1}.
$

Conversely if there exists $C_{\p}>0$ such that 
$Cap_{\o}(\p <-t) \leq C_{\p} /t^{p+2}$ for all $t \geq 1$, 
then $\f \in {\mathcal E}^q(X,\o)$ for all $q<p$.
\end{lem}

\begin{proof}
It follows from the comparison principle that for all $t \geq 1$,
\begin{equation}
\int_{(\f<-t)} \o_{\f}^2 \leq t^2 Cap_{\o}(\f<-t).
\end{equation}
Indeed set $\f_s:=\max(\f,-s)$, where $s>t$. Then
$\f_s \in PSH(X,\o) \cap L^{\infty}(X)$ and $\f_s=\f$ near the boundary
of $(\f<-t)$, hence
$$
\int_{(\f<-t)} \o_{\f}^2=\int_{(\f<-t)} \o_{\f_s}^2.
$$
Consider $u=\f_s/s \in PSH(X,\o)$. Then $-1 \leq u \leq 0$ and 
$s^{-2} \o_{\f_s}^2 \leq \o_u^2$, hence
$$
\frac{1}{s^2} \int_{(\f<-t)} \o_{\f_s}^2 \leq \int_{(\f<-t)} \o_u^2
\leq Cap_{\o}(\f<-t),
$$
which yields (6) by letting $s \searrow t$. We infer
\begin{eqnarray*}
0 &\leq& \int_X (-\f)^p \o_{\f}^2 = p\int_1^{+\infty} t^{p-1} \o_{\f}^2(\f<-t) dt \\
&\leq& p \int_1^{+\infty} t^{p+1} Cap_{\o}(\f<-t) dt 
= \frac{p}{p+2} \int_X (-\f)^{p+2} dCap_{\o}
\end{eqnarray*}

The second inequality also follows from the comparison principle. We need to estimate
$Cap_{\o}(\f<-t)$ from above, for $t \geq 1$. Fix $u \in PSH(X,\o)$
with $-1 \leq u \leq 0$ and observe that
$\f/t \in PSH(X,\o)$ with
$$
\left( \f<-2t \right) \subset \left( \f/t< u-1 \right)
\subset \left( \f <-t \right).
$$
Therefore $\int_{(\f<-2t)} \o_u^2 \leq \int_{(\f<-t)} \o_{\f/t}^2$. 
Now $\o_{\f/t} \leq \o+t^{-1} \o_{\f}$, thus
\begin{equation}
Cap_{\o}(\f<-2t) \leq \left[ \o^2+\frac{2}{t} \o \wedge \o_{\f}
+\frac{1}{t^2} \o_{\f}^2 \right] (\f<-t).
\end{equation}
This yields the second inequality. The remaining assertions are straightforward
consequences of (6), (7) and Chebyshev inequality.
\end{proof}

These estimates allow us to give now several examples of functions which belong
to the classes ${\mathcal E}^p(X,\o)$.

\begin{exa}
\text{ }

1) Assume $X=\P^2$ is the complex projective space and $\o=\o_{FS}$ is the Fubini-Study
K\"ahler form. We let $[z_0:z_1:z_2]$ denote the homogeneous coordinaets on $\P^2$. Consider 
$\f[z_0:z_1:z_2]=\log ||(z_1,z_2)||-\log ||(z_0,z_1,z_2)|| \in {\mathcal E}(X,\o)$. This is
a ${\mathcal C}^{\infty}$-smooth function in $\P^2 \setminus \{a\}$, where
$a=[1:0:0]$, which has Lelong number $1$ at point $a$, hence 
$\f \notin {\mathcal E}^1(X,\o)$. One can compute explicitly
$$
Cap_{\o}(\f<-t) \simeq \frac{C}{t^2}
$$
by comparing the capacity $Cap_{\o}$ with the local Monge-Amp\`ere capacity of
Bedford and Taylor near point $a$ (see [16]).

Consider now $\f_{\a}:=-(-\f)^{\a}$ for $0 <\a<1$. Then 
$\f_{\a} \in {\mathcal E}^p(X,\o)$ for $\a<2/(p+2)$, as follows from lemma 6.2.

2) This first example can be slightly generalized as follows. Let
$\f$ be any $\o$-psh function such that $\f \leq -1$. Then 
$\f_{\a}:=-(-\f)^{\a} \in PSH(X,\o)$ whenever $0 \leq \a \leq 1$. 
It follows moreover from lemma 5.2 that $\f_{\a} \in {\mathcal E}^p(X,\o)$
as soon as $\a<1/(p+2)$.
In particular every pluripolar set is included in the $-\infty$ locus
of a function $\f \in {\mathcal E}^p(X,\o)$, for all $p \geq 1$.

In the same vein observe that $\Phi:=-\log(-\f) \in PSH(X,\o)$ if
$\f\leq -2$. Moreover $\Phi \in {\mathcal E}^p(X,\o)$ for all $p \geq 1$,
although $\Phi$ is not necessarily bounded.

3) Assume $X=\P_x^1 \times \P_y^1$ and $\o(x,y)=\a(x)+\a(y)$, where 
$\a$ denotes the Fubini-Study K\"ahler form on $\P^1$.
Assume $\f(x,y)=u(x)+v(y)$ where $u,v$ come both 
from $\a$-psh functions on $\P^1$, with $u$ smooth while $v$ is singular.
Then $\o_{\f}=\a_u(x)+\a_v(y)$ and $\o_{\f}^2=2\a_u(x) \wedge \a_v(y)$ so that
$$
\f \in L^p(\o_{\f}^2) \Leftrightarrow v \in L^p(\a_v)
\Leftrightarrow \f \in L^p(\o \wedge \o_{\f}).
$$
In particular one can get $\f \in L^p(\o_{\f}^2)$ for some $p \geq 1$
but $\f \notin L^{p+\e}(\o \wedge \o_{\f})$ whenever $\e>0$.
\end{exa}

Observe that there are functions in examples 1 and 3 above that belong to 
${\mathcal E}^p(X,\o)$ for some $p \geq 1$, but not to ${\mathcal E}^{p+1}(X,\o)$.
The last example explains partially why there is a gap in the estimates given
by lemma 5.2: we have the right exponent when $\f \in {\mathcal E}^p(X,\o)$
also satisfies $\f \in L^{p+1}(\o \wedge \o_{\f})$, but this integrability
condition is not necessarily satisfied unless 
$\f \in {\mathcal E}^{p+1}(X,\o)$.
It is satisfied  however, when the function $\f$ has singularities
in a ``small compact'' (see example 1.4.2).
More precisely we have the following:

\begin{pro}
Assume $\f \in {\mathcal E}^p(X,\o)$ is bounded near some ample divisor.
Then $\f \in L^{p+1}(\o \wedge \o_{\f})$.
\end{pro}

\begin{proof}
Let $\f \in {\mathcal E}^p(X,\o)$, $\f \leq 0$, be bounded in a neighborhood
$V$ of some ample divisor $D$. We assume for simplicity that the current
$[D]$ of integration along $D$ is cohomologous to $\o$. Let $\o' \geq 0$ be a smooth
closed $(1,1)$-form cohomologous to $\o$ such that $\o' \equiv 0$ in $X \setminus V$.
Fix $\chi \leq 0$ such that $\o=\o'+dd^c \chi$. 
Assume first $\f$ is smooth. It follows from Stokes theorem that
\begin{eqnarray*}
0 \leq \int_X (-\f)^{p+1} \o \wedge \o_{\f} &=& 
\int_X (-\f)^{p+1} \o'\wedge \o_{\f} +\int_X (-\f)^{p+1} \o_{\f} \wedge dd^c \chi \\
\!\!\!\!\!\!\!\!
&\leq& \!\!\!\!
||\f||_{L^{\infty}(V)}^{p+1}+(p+1) ||\chi||_{L^{\infty}(X)} \int_X (-\f)^p \o_{\f}^2
<+\infty.
\end{eqnarray*}
Now we can get a similar control on $||\f||_{L^{p+1}(\o \wedge \o_{\f})}$ by approximating
$\f$ by a decreasing sequence of smooth $\o$-psh functions and by
using theorem 4.4.
\end{proof}

We now want to give some examples of probability measures which can be expressed
as $\o_{\p}^2$, $\p \in {\mathcal E}^p(X,\o)$.
Observe first that if $PSH(X,\o) \subset L^1(\mu)$, then in particular
${\mathcal E}^1(X,\o) \subset L^1(\mu)$, hence $\mu=\o_{\p}^2$ for some
function $\p \in {\mathcal E}^1(X,\o)$. Every measure which decomposes as
$$
\mu=\Theta+dd^c (R),
$$
where $\Theta$ is a smooth form and $R$ is a {\it positive} current of bidimension 
$(1,1)$, satisfies $PSH(X,\o) \subset L^1(\mu)$. There are several examples of
such measures arising in complex dynamics [15]. We don't know if
$\p$ is necessarily bounded in this case (this is trivially true in dimension $1$).

When $\mu=f\o^2$ has density $f \in L^r(X)$, $r>1$, S.Kolodziej has proved [18]
that $\mu=\o_{\p}^2$ for some {\it bounded} $\o$-psh function $\p$.
This is because $\mu$ is strongly dominated by $Cap_{\o}$ 
in this case (see proposition 6.5 below).
When the density is only in $L^1$, this does not work. Consider for instance
$\mu=f\o^2$, where $f \in {\mathcal C}^{\infty}(X \setminus \{a\})$ is
such that 
$$
f(z) \simeq \frac{1}{||z||^4 (-\log ||z||)^2}-1
$$
near the point $a$($=0$ in a local chart). Observe that
$$
\f(z):=\e \chi(z) \log ||z|| \in PSH(X,\o)
$$
if $\chi$ is a cut-of function such that $\chi \equiv 1$ near $a=0$, and $\e>0$ is small 
enough. Now $\f \notin L^1(\mu)$ but still
${\mathcal E}^1(X,\o) \subset L^1(\mu)$.
Therefore there exists $\p \in {\mathcal E}^1(X,\o)$ such that $\mu=\o_{\p}^2$, as
follows from theorem 5.1. Note however that $\p \notin {\mathcal E}^p(X,\o)$
for $p>2$. Indeed $\f_{\a}:=-(-\f)^{\a} \in {\mathcal E}^p(X,\o)$ if
$\a<2/(p+2)$ (see example 6.3.1) and $\f_{\a} \in L^p(\mu)$ implies
$\a p<1$, hence $p \leq 2$.

Observe also that there are measures $\mu=f \o^2$ with $L^1$-density such that
${\mathcal E}^1(X,\o) \not\subset L^1(\mu)$: one can consider for instance
$f_{\e}$ that looks locally near $a=0$ like
$[ ||z||^4 (-\log ||z||)^{1+\e}]^{-1}$, for $\e>0$ small enough.

In order to give further examples, we need to relate integrability properties 
of $\mu$ to the way it is dominated by $Cap_{\o}$. This is the contents of the following:

\begin{pro}
Let $\mu$ be a probability measure on $X$.

Assume there exists $\a>p/(p+1)$ and $A>0$ such that
\begin{equation}
\mu(E) \leq A Cap_{\o}(E)^{\a},
\end{equation}
for all Borel set $E \subset X$. Then ${\mathcal E}^p(X,\o) \subset L^p(\mu)$.

Conversely assume ${\mathcal E}^p(X,\o) \subset L^p(\mu)$. Then there exists
$0<\a<1$ and $A>0$ such that $(8)$ is satisfied.
\end{pro}

\begin{proof}
We can assume w.l.o.g. that $\a \leq 1$. Let $\f \in {\mathcal E}^p(X,\o)$
with $\sup_X \f=-1$. It follows from H\"older inequality that
\begin{eqnarray*}
\lefteqn{ \! \! \! \! \! \! \! \! \! 
0  \leq  \int_X (-\f)^p d\mu = p \int_1^{+\infty} t^{p-1} \mu(\f<-t) dt} \\
&\leq& pA \int_1^{+\infty}t^{p-1} \left[ Cap_{\o}(\f<-t) \right]^{\a} dt \\
&\leq& pA \left[ \int_1^{+\infty} t^{\frac{p-\a(p+1)}{1-\a}-1} dt \right]^{1-\a}
\cdot \left[ \int_1^{+\infty} t^p Cap_{\o}(\f<-t) dt \right]^{\a}.
\end{eqnarray*}
The first integral in the last line converges since $p-\a(p+1)<0$ and the last is dominated
by $C_{p,\a} \left( \int (-\f)^p \o_{\f}^2 \right)^{\a}$ by lemma 6.2.
Therefore ${\mathcal E}^p(X,\o) \subset L^p(\mu)$.

Assume conversely that ${\mathcal E}^p(X,\o) \subset L^p(\mu)$.
It follows from Theorem 5.1 that $\mu=\o_{\p}^2$ for a unique
function $\p \in {\mathcal E}^p(X,\o)$ such that $\sup_X \p=-1$.
We claim then that there exists $\g_p \in ]0,1[$ and $A>0$ such that
for all functions $\f \in PSH(X,\o)$ with $-1 \leq \f \leq 0$, one has
\begin{equation}
0 \leq \int_X (-\f)^p \o_{\p}^2 \leq A \left[ \int_X (-\f)^p \o_{\f}^2 \right]^{\g_p}.
\end{equation}
We leave the proof of this claim to the reader: an application of 
H\"older's inequality in the style of section 4 yields $\g_p=(1-1/p)^2$
if $p>1$, while Cauchy-Schwarz inequality (in the style of lemma 3.6)
yields $\g_1=1/4$. We apply now (9) to the extremal function
$\f=h_{E,\o}^*$ introduced in [16]. It follows from
theorem 3.2 in [16] that
$$
0 \leq \mu(E) \leq \int_X (-h_{E,\o}^*)^p d\mu \leq A  Cap_{\o}(E)^{\g_p}.
$$
\end{proof}

This proposition allows to produce several examples of measures satisfying
${\mathcal E}^p(X,\o) \subset L^p(\mu)$ as in the local theory
(see [18], [26]).
It can also be used, together with theorem 5.1, to prove that functions from the 
local classes of Cegrell ${\mathcal E}^r(\Omega)$, $\Omega$ a bounded hyperconvex
domain of $\C^2$, can be sub-extended as global functions
$\f \in {\mathcal E}^{p_r}(\P^2,\o_{FS})$
(see [10] for similar results).

\section{Generalizations and applications}

In this section we assume $X$ is a compact K\"ahler manifold
of arbitrary dimension $k \geq 1$, equipped with a K\"ahler form $\o$
such that $\int_X \o^k=1$.

\subsection{Higher dimension}

Our main results easily extend to higher dimension.
Apart from the slightly more involved computations, one difficulty
we face is that of the definition of the complex Monge-Amp\`ere operator
(the $L^2$- condition on the gradient is not sufficient in dimension
$\geq 3$). We adopt here an ad hoc definition.

\begin{defi}
Fix $p \geq 1$.
We let ${\mathcal E}^p(X,\o)$ denote the set of functions $\f \in {\mathcal E}(X,\o)$
such that there exists a sequence $\f_j \in PSH(X,\o) \cap L^{\infty}(X)$ with
$$
\f_j \searrow \f \; \text{ and } \; 
\sup_{j \in \N} \left( \int_X |\f_j|^p \o_{\f_j}^k \right) < +\infty.
$$
\end{defi}

One can then prove similar results that those obtained in section 4,
where $X$ had dimension 2. Note that most constants will then depend
on $k=\dim_{\C} X$, since one has to integrate by parts several times.
As an example lemma 4.2 becomes:

\begin{lem}
Let $\f,\p \in PSH(X,\o) \cap L^{\infty}(X)$ 
with $\f \leq \p \leq 0$. 
Fix $p \geq 1$ and $0 \leq j \leq k$.
Then

1) $0 \leq \int_X (-\f)^p \o^k \leq \int_X (-\f)^p \o^{k-j} \wedge \o_{\f}^j
\leq \int_X (-\f)^p \o_{\f}^k$.

2) $0 \leq \int_X (-\p)^p \o^{k-j} \wedge \o_{\p}^j \leq (p+1)^j 
\int_X (-\f)^p \o^{k-j} \wedge \o_{\f}^j$.

\noindent In particular
$$
0 \leq \int_X (-\p)^p \o_{\p}^k \leq (p+1)^k \int_X (-\f)^p \o_{\f}^k.
$$
\end{lem}

The proof is essentially the same as that of lemma 4.2. 
The result corresponding to Theorem 4.4 then shows that the
complex Monge-Amp\`ere operator is well-defined in the class
${\mathcal E}^p(X,\o)$ and is continuous on decreasing sequences.

\begin{thm}
Let $\f \in {\mathcal E}^p(X,\o)$.
There exists a probability measure $\mu_{\f}$
such that if  $(\f_j)$ is any sequence of 
bounded $\o$-psh functions that decreases towards $\f$,
then $\sup (\int |\f_j|^p \o_{\f_j}^2)<+\infty$ and 
for all $0 \leq q<p$ ,
$$
(-\f_j)^q \o_{\f_j}^k \longrightarrow (-\f)^q \mu_{\f}.
$$

The convergence still holds if $(\f_j)$ merely converges
to $\f$ in capacity.
\end{thm}

Of course we denote by $\o_{\f}^k$ the measure $\mu_{\f}$.
The proof is a straighforward generalization of that of Theorem 4.4:
observe that the monotonicity implies convergence in capacity.
There is again unicity of solutions to Monge-Amp\`ere equations in
the class ${\mathcal E}^1(X,\o)$:

\begin{thm}
Let $\f,\p \in {\mathcal E}^1(X,\o)$ be such that $\o_{\f}^k \equiv \o_{\p}^k$.
Then $\f-\p$ is constant.
\end{thm}

Since the proof is technically more involved than that
of Theorem 3.4, but relies on the same ideas, we only sketch it.
\vskip.1cm

\begin{sketch}
The idea is again to prove that $\nabla f=0$, where $f:=\f-\p$.
We want to get a bound from above $\int df \wedge d^c f \wedge \o^{k-1}$
by a function of $\int (-\f) \o_{\f}^k, \int (-\p) \o_{\p}^k$, and
$\int df \wedge d^c f \wedge \sum_{\a+\b=k-1} \o_{\f}^{\a} \wedge \o_{\p}^{\b}$,
which vanishes if this last integral is zero. We will be done
since by Stokes Theorem,
$$
\int df \wedge d^c f \wedge 
\left( \sum_{\a+\b=k-1} \o_{\f}^{\a} \wedge \o_{\p}^{\b} \right)
=\int -f (\o_{\f}^k -\o_{\p}^k)=0.
$$

Let $T$ be a positive closed current of bidegree $(k-2,k-2)$. The clue of the
proof lies, as in 3.4, in the following integrations by parts.
Observe that
$$
\int df \wedge d^c f \wedge \o \wedge T=
\int df \wedge d^c f \wedge \o_{\f} \wedge T
-\int df \wedge d^c f \wedge dd^c \f \wedge T=I+II.
$$
Now 
$$
II=\int df \wedge d^c \f \wedge \o_{\f} \wedge T-\int df \wedge d^c \f \wedge \o_{\p} \wedge T
=II'+II'',
$$
and we estimate II' and II'' by using Cauchy-Schwarz inequality,
$$
II' \leq \left( \int df \wedge d^c f \wedge \o_{\f} \wedge T \right)^{1/2}
\cdot 
\left( \int d\f \wedge d^c \f \wedge \o_{\f} \wedge T \right)^{1/2},
$$
and similarly for $II''$. We now use these estimates by putting
$T=\o^{k-2}$, then $T=\o^{k-3} \wedge \o_{\f}$ or $T=\o^{k-3}  \wedge \o_{\p}$,
\ldots, $T=\o^a \wedge \o_{\f}^b \wedge \o_{\p}^c$.
Each term of the form $\int d\f \wedge d^c \f \wedge \o_{u} \wedge T$, $u=\f$ or $\p$, 
will be dominated from above
by $\max ( \int (-\f) \o_{\f}^k; \, \int (-\p) \o_{\p}^k)$, thanks to the k-dimensional
version of proposition 3.2. Each term of the form
$\int df \wedge d^c f \wedge \o_{u}  \wedge T$ will eventualy be dominated
by integrals of the form
$\int df \wedge d^cf \wedge \o_{\f}^j \wedge \o_{\p}^{k-1-j}$ which are all
bounded from above by
$\int df \wedge d^c f \wedge \sum_{\a+\b=k-1} \o_{\f}^{\a} \wedge \o_{\p}^{\b}$.
\end{sketch}

All other results from previous sections generalize easily and yield
the following higher dimensional version of Theorem 5.1:

\begin{thm}
Let $\mu$ be a probability measure on $X$ and fix $p \geq 1$.
Then there exists a unique function $\p \in {\mathcal E}^p(X,\o)$ such that
$$
\mu=\o_{\p}^k \; \; \text{ and } \; \; \sup_X \p=-1
$$
if and only if ${\mathcal E}^p(X,\o) \subset L^p(\mu)$.
\end{thm}

\subsection{Complex dynamics}

Functions of the class ${\mathcal E}^1(X,\o)$ and probability measures
$\mu$ such that ${\mathcal E}^1(X,\o) \subset L^1(\mu)$ naturally arise 
in complex dynamics. Namely let 
$f:X \rightarrow X$ be a meromorphic endomorphism whose topological degree
$d_t(f)$ is large in the sense that
$d_t(f)>\l_{k-1}(f)$, where $\l_{k-1}(f)$ denotes the spectral radius
of the linear action induced by $f^*$ on $H^{k-1,k-1}(X,\R)$.
In this case there exists a unique invariant measure $\mu_f$
of maximal entropy which can be decomposed as
$$
\mu_f:=\Theta+dd^c ({\mathcal T}),
$$
as was proved by the first author in [15]. Here $\Theta$
is a smooth probability measure and ${\mathcal T} \geq 0$ is
a positive current of bidegree $(k-1,k-1)$. In particular
$$
{\mathcal E}^1(X,\o) \subset PSH(X,\o) \subset L^1(\mu) ,
$$
as follows from Stokes theorem (see Theorem 2.1 in [15]
and Example 2.8 in [16]).
It follows therefore from our main result (Theorem 5.1) that there
exists a unique function $g_f \in {\mathcal E}^1(X,\o)$ such that
$\sup_X g_f=-1$ and 
$$
\mu_f=(\o+dd^c g_f)^k.
$$
It is an interesting problem to establish further regularity properties
of $g_f$ in order e.g. to estimate the pointwise dimension of the measure $\mu_f$.

In a similar but slightly different direction, let us consider now the case
where $X=\P^2$ is the complex projective plane 
equipped with the Fubini-Study K\"ahler form $\o$, 
and $f:\P^2 \rightarrow \P^2$ is {\it birational}, i.e. $d_t(f)=1$, 
with $\l_1(f)>1$. When $f$ satisfies a technical -- but generic -- condition, 
E.Bedford and J.Diller have constructed in [2] a canonical invariant probability
measure
$$
\mu_f:=\o_{g^+} \wedge \o_{g^-}, 
$$
where the dynamical Green functions $g^{\pm} \in {\mathcal E}(\P^2,\o)$
are such that $g^{\pm} \in L^1(\mu_f)$. The functions
$g^+,g^-$ do not belong to the class ${\mathcal E}^1(\P^2,\o)$ because
they have positive Lelong numbers at points of indeterminacy of the mappings 
$f^n$, $n \in \Z$, however we have the following result.

\begin{pro}
Set $g:=\max (g^+,g^-)$. Then
$g \in {\mathcal E}^1(\P^2,\o)$.
\end{pro}

\begin{proof}
We can assume without loss of generality that $g^+,g^- \leq 0$, hence $g \leq 0$.
It follows from Proposition 1.3.3 that $g \in {\mathcal E}(X,\o)$,
hence the Monge-Amp\`ere measure $\o_g^2$ is well defined. Now
$$
0 \leq \int_X (-g) \o_g^2 \leq \int_X (-g^+) \o_g^2=\int (-g^+) \o\wedge \o_g
+\int_X (-g^+) dd^c g \wedge \o_g.
$$
The first integral in the RHS is finite thanks to Proposition 1.3.1.
The last one can be bounded from above, using Stokes theorem and
Proposition 1.3.1 again,
$$
\int_X (-g^+) dd^c g \wedge \o_g=\int_X (-g) dd^c g^+ \wedge \o_g
\leq \int_X (-g) \o_{g^+} \wedge \o_g+O(1)
$$
Using a similar integration by parts, we obtain
$$
\int_X (-g) \o_{g^+} \wedge \o_g \leq \int_X (-g^-) \o_{g^+} \wedge \o_g
\leq \int_X (-g) \o_{g^+} \wedge \o_{g^-}+O(1),
$$
hence
$$
\int_X (-g) \o_g^2 \leq \int_X (-g) \o_{g^+} \wedge \o_{g^-}+O(1)
\leq \int_X (-g^+)d\mu_f+O(1) < +\infty.
$$
\end{proof}

It is an interesting problem to determine whether 
$\mu_f=(\o+dd^cg)^2$. This is the case when e.g. $f$ is a complex H\'enon mapping,
and it would imply -- by the converse to our main result Theorem 5.1 -- that
${\mathcal E}^1(\P^2,\o) \subset L^1(\mu_f)$, hence 
in particular $\mu_f$ does not charge
pluripolar sets (by Example 6.3.2).

\subsection{K\"ahler-Einstein metrics}

It is well-known that solving Monge-Amp\`ere equations
$$
\text{ }
\!\!\!\!\!\!  \!\!\! \!\!\! \!\!\! \!\!\! \!\!\! 
[MA](X,\o,\mu) \hskip2cm
(\o+dd^c \f)^k=\mu    
$$
is a way to produce K\"ahler-Einstein metrics (see [8], [1], [25], [22], [18]).
In the classical case, the measure $\mu=f \o^k$ admits a smooth density
$f>0$. When the ambient manifold has some singularities (which is often the case
in dimension $\geq 3$), one has to allow the equation
$[MA](X,\o,\mu)$ to degenerate in two different ways:
resolving the singularities $\pi:\tilde{X} \rightarrow X$ of $X$ yields
a new equation $[MA](\tilde{X},\tilde{\o},\tilde{\mu})$, where
\begin{enumerate}
\item $\{\tilde{\o}\}=\{ \pi^* \o\}$ is a semi-positive and big class (one looses
strict positivity along the exceptional divisors);
\item $\tilde{\mu}=\tilde{f} \tilde{\o}^k$ is a measure with density 
$0 \leq \tilde{f} \in L^p$, $p>1$, which may have zeroes and poles along some of 
the exceptional divisors.
\end{enumerate}

We have focused in this paper on the second type of degeneracy. We would like to 
mention that our techniques are supple enough so that we can produce solutions
$\f \in{\mathcal E}^1(\tilde{X},\tilde{\o})$ to the Monge-Amp\`ere equations
$[MA](\tilde{X},\tilde{\o},\tilde{\mu})$, even when $\{\tilde{\o}\}$ is
merely big and semi-positive rather than K\"ahler. We will develop
this in our forthcoming article [14].

\vskip .2cm

Vincent Guedj \& Ahmed Zeriahi

Laboratoire Emile Picard

UMR 5580, Universit\'e Paul Sabatier

118 route de Narbonne

31062 TOULOUSE Cedex 04 (FRANCE)

guedj@picard.ups-tlse.fr

zeriahi@picard.ups-tlse.fr

\end{document}